\documentclass[12pt, a4paper]{article}
\usepackage[pagewise]{lineno}
\usepackage{graphicx}
\usepackage{amssymb}
\usepackage{latexsym, bm}
\usepackage{multicol}
\usepackage{indentfirst}
\usepackage{amssymb,amsfonts}
\usepackage{amsmath}
\usepackage{cases}
\usepackage[noadjust]{cite}
\usepackage[colorlinks,citecolor=red]{hyperref}
\newtheorem{theorem}{Theorem}
\newtheorem{lemma}{Lemma}

\newtheorem{proposition}{Proposition}
\newtheorem{definition}{Definition}

\newtheorem{remark}{Remark}

\usepackage{amssymb}
\numberwithin{equation}{section}
\numberwithin{theorem}{section}
\numberwithin{remark}{section}
\numberwithin{definition}{section}
\numberwithin{lemma}{section}
\numberwithin{corollary}{section}
\numberwithin{proposition}{section}
\DeclareMathOperator{\dive}{div}
\title{On the energy equality for the incompressible viscoelastic flows}
\author{Wenke Tan\footnote{tanwenkeybfq@163.com}\quad Fan Wu\footnote{wufan0319@yeah.net;wufan0319@gmail.com}
	\\
{\small Key Laboratory of Computing and Stochastic Mathematics (Ministry of Education),}\\
{\small School of Mathematics and Statistics, Hunan Normal University,}\\
{\small Changsha, Hunan 410081, China}\\
}

\date{}

\begin{document}
\maketitle
{\bf Abstract:}
In this paper, we study the problem of energy conservation for the solutions to the incompressible viscoelastic flows. First, we consider Leray-Hopf weak solutions in the bounded Lipschitz domain $\Omega$ in $\mathbb{R}^d\,\, (d\geq 2)$. We prove that under the Shinbrot type conditions $
u \in L^{q}_{loc}\left(0, T ; L^{p}(\Omega)\right)
\text { for any } \frac{1}{q}+\frac{1}{p} \leq \frac{1}{2}, \text { with } p \geq 4,\text{ and }  {\bf F} \in L^{r}_{loc}\left(0, T ; L^{s}(\Omega)\right) \text { for any } \frac{1}{r}+\frac{1}{s} \leq \frac{1}{2}, \text { with } s \geq 4  $, the boundary conditions $u|_{\partial\Omega}=0,\,\,{\bf F}\cdot n|_{\partial\Omega}=0$ can inhibit the boundary effect and guarantee the validity of energy equality. Next, we apply this idea to deal with the case $\Omega= \mathbb{R}^d\,\,(d=2, 3, 4)$, and showed that the energy is conserved for 	
$u\in L_{loc}^{q}\left(0,T;L_{loc}^{p}\left(\mathbb{R}^{d}\right)\right)$ with $ \frac{2}{q}+\frac{2}{p}\leq1, p\geq 4
$
and
$
	{\bf F}\in L_{loc}^{r}\left(0,T;L_{loc}^{s}\left(\mathbb{R}^{d}\right)\right)\cap L^{\frac{4d+8}{d+4}}\left(0,T;L^{\frac{4d+8}{d+4}}\left(\mathbb{R}^{d}\right)\right)$  with $\frac{2}{r}+\frac{2}{s}\leq1,  s\geq 4
$. This result shows that the behavior of solutions in the finite regions and the behavior at infinite play different roles in the energy conservation. Finally, we consider the problem of energy conservation for distributional solutions and show energy equality for the distributional solutions belonging to the so-called Lions class $L^4L^4$.

\medskip

\section{Introduction}
The Oldroyd-type models capture the rheological phenomena of both the fluid motions and the elastic features of non-Newtonian fluids. We study the simplest case in which the relaxation and retardation times are both infinite. More specifically, we consider the following system of equations for an incompressible viscoelastic fluid in $\Omega\subseteq\mathbb{R}^d\, (d\geq 2)$:
\begin{equation}\label{1.1}
\left\{
             \begin{array}{lr}
             \partial_t u+(u\cdot\nabla)u-\mu\Delta u+\nabla P=\nabla\cdot({\bf F}{\bf F}^\top),& \\
            \partial_t {\bf F}+(u\cdot\nabla){\bf F}=(\nabla u){\bf F},&\\
             \nabla\cdot u=0,&\\
             u|_{\partial\Omega}=0,\,\, n_j{\bf F}_{jk}|_{\partial\Omega}=0,\\
             u(x,0)=u_{0}(x),{\bf F}(x,0)={\bf F}_{0}(x),&
\end{array}
\right.
\end{equation}
where $u$ is the velocity field, ${\bf F}$ is the local deformation tensor of the fluid, and $P$ represents the pressure. The constant $\mu>0$ is the kinetic viscosity. Here $(\nabla\cdot({\bf F}{\bf F}^\top))_i=\partial_j({\bf F}_{ik}{\bf F}_{jk})$ and $(\nabla u)_{ij}=\partial_ju_i$. For convenience, we assume $\mu=1$ throughout this paper. The system \eqref{1.1} is also called as
the viscoelastic Navier-Stokes equations.

The complicated rheological phenomena of complex fluids is a consequence of interactions
between the (microscopic) elastic properties and the (macroscopic) fluid motions \cite{RGL,RM}. The system \eqref{1.1} serves as an important model for the study of the dynamics of complex fluids. It
presents the competition between the elastic energy and the kinetic energy, while the deformation tensor ${\bf F}$ carries all the transport/kinematic information of the micro-structures and configurations in complex fluids \cite{LLZ,LZ,LXZ}. It is worth pointing out that $\nabla\cdot {\bf F^{\top}}$ satisfies the following system
\begin{align*}
\partial_t(\nabla\cdot {\bf F^{\top}})+u\cdot\nabla(\nabla\cdot {\bf F^\top})=0,
\end{align*}
this means $\nabla\cdot {\bf F^\top}=0$ if one imposes the initial compatible condition $\nabla\cdot {\bf F}_{0}^\top(x)=0$. This information plays  a
key role in proving global existence for smooth small data. However, if we are only concerned with the energy conservation problem, we shall show that the compatibility condition $\nabla\cdot {\bf F^\top}=0$ is not necessary for energy equality.

In the last few decades, there are many research results on the system \eqref{1.1}.
Global existence for solutions near some equilibrium of the viscous analog of the system \eqref{1.1} has been verified by Lin-Liu-Zhang \cite{LLZ} for two-dimensional flow and Lei-Liu-Zhou \cite{LZL} for three-dimensional flow. Hu-Wu \cite{HW} provide an alternative proof of the existence of global smooth solutions near some equilibrium and obtained the optimal $L^2$ decay rates of smooth solutions. Moreover, the weak-strong uniqueness property in the class of finite energy weak solutions was established. Recently, Hu and Lin in \cite{HL} proved the global-in-time existence of the 2D Leray-Hopf type weak solutions in the physical energy space.
For a further discussion of these topics, we may refer to the survey paper written by Lin \cite{LINFANGHUA}. In addition,  concerning the regularity criteria of the system \eqref{1.1},
see \cite{HH,FZZ,YBQ} for some related discussions. Besides, it is also worth mentioning that there is some recent progress on the compressible viscoelastic flows, see \cite{HXP,HUW,CW,IY} and the references therein.

System \eqref{1.1} can be regarded as a coupling of a parabolic system with a hyperbolic one. For the smooth solution, it (formally) obeys the following basic energy law:
\begin{equation}\label{1.2}
\frac{1}{2}\int_{\Omega}\left(|u|^{2}+|{\bf F}|^{2}\right)(t) d x+\int_{0}^{t} \int_{\Omega}|\nabla u|^{2} d x d s=\frac{1}{2} \int_{\Omega}\left(\left|u_{0}\right|^{2}+\left|{\bf F}_{0}\right|^{2}\right) d x.
\end{equation}
However, for weak solutions with less regularity, identity \eqref{1.2} may be fail. That is, a weak solution $(u,\bf F)$ to \eqref{1.1} satisfies the following
energy inequality:
\begin{equation}\label{EI}
	\frac{1}{2}\int_{\Omega}\left(|u|^{2}+|{\bf F}|^{2}\right)(t) d x+\int_{0}^{t} \int_{\Omega}|\nabla u|^{2} d x d s\leq\frac{1}{2} \int_{\Omega}\left(\left|u_{0}\right|^{2}+\left|{\bf F}_{0}\right|^{2}\right) d x.
\end{equation}
From the physical point of view we would expect that a weak solution $(u, \bf F)$ satisfies
the energy equality \eqref{1.2}.
It is natural to ask the following questions:

$\bullet$ Which regularity of Leray-Hopf weak solutions for the validity of the energy equality \eqref{1.2}?

$\bullet$ Which regularity of distributional solutions (see definition \ref{de2}) possible to guarantee that energy equality \eqref{1.2} holds?

For the incompressible viscoelastic fluids problem \eqref{1.1}, if let $\mu=0$, then the equations reduce to the ideal viscoelastic flow. The famous Onsager conjecture for the Euler equations predicts the threshold regularity for energy conservation. In this direction,  Onsager \cite{ON} considers periodic 3-dimensional weak solutions of the incompressible Euler equations, where the velocity $u$ satisfies the uniform H\"older condition
$$
|u(x, t)-u(x', t)|\leq C|x-x'|^\alpha,
$$
for constants $C$ and $\alpha$ independent of $x, x'$ and $t$.

$\bullet$  If $ \alpha > \frac{1}{3}$
, then the total kinetic energy $E(t) = \frac{1}{2}\int |u(x, t)|^2 dx$ is a constant.

 $\bullet$  For any $ \alpha < \frac{1}{3}$
, there is a solution $u$ for which the kinetic energy $\frac{1}{2}\int|u(x,t)|^2dx$ is not a constant.\\
For $\alpha > \frac{1}{3}$, Eyink \cite{EGL} and Constantin-E-Titi \cite{CPT} proved the energy conservation. Recently,  Isett \cite{IS} gave the complete answer to the second part of the conjecture.
Cheskidov-Constantin-Friedlander-Shvydkoy in \cite{CCFS} proved that weak solution in the following class
\begin{equation}\label{1.3}
	u\in L^{3}\left(0, T; B_{3, c(\mathbb{N})}^{1 / 3}\right) \cap C_{w}\left(0, T; L^{2}\right)
\end{equation}
conserves energy. Bardos-Titi \cite{BT} considered the Onsager's conjecture for the incompressible
Euler equations in bounded domains and showed that, in a bounded domain $\Omega \subset \mathbb{R}^{d}$, with $\partial \Omega \in C^{2}$, any weak solution $(u(x, t), p(x, t))$, of the Euler equations of ideal incompressible fluid in $\Omega \times(0, T) \subset \mathbb{R}^{d} \times \mathbb{R}_{t}$, with the impermeability boundary condition $u \cdot \vec{n}=0$ on $\partial \Omega \times(0, T)$, is of constant energy on the interval $(0, T)$, provided the velocity field $u \in L^{3}\left((0, T) ; C^{0, \alpha}(\bar{\Omega})\right)$, with $\alpha>\frac{1}{3}$. Concerning the energy conservation of ideal viscoelastic flows, two sufficient conditions of the energy conservation for weak solutions of incompressible viscoelastic flows are established in \cite{HZ}. In particular, for a periodic domain in $\mathbb{R}^{3}$, energy conservation is proved for $u$ and ${\bf F}$ in certain Besov spaces. Furthermore, in the whole space $\mathbb{R}^{3}$, they showed that the conditions on the velocity $u$ and the deformation tensor ${\bf F}$ can be relaxed, that is, $u \in B_{3, c(\mathbb{N})}^{\frac{1}{3}}$ and ${\bf F} \in B_{3, \infty}^{\frac{1}{3}}$. In fact, for the system \eqref{1.1} ($\mu>0$), They obtained that the energy is conserved for $u \in L^{r}\left(0, T ; L^{s}(\mathbb{T}^d)\right)$ for any $\frac{1}{r}+\frac{1}{s} \leqslant \frac{1}{2}$, with $s \geqslant 4$, and ${\bf F} \in L^{m}\left(0, T ; L^{n}(\mathbb{T}^d)\right)$ for any $\frac{1}{m}+\frac{1}{n} \leqslant \frac{1}{2}$, with $n \geqslant 4$. However, they failed to prove the energy equality of \eqref{1.1} in a bounded domain $\Omega\subset \mathbb{R}^d$ because they could not deal with the boundary effect.

Notice that in the absence of deformation tensor (i.e., ${\bf F}=0$) and $\mu>0$, system \eqref{1.1} reduces to the incompressible Navier-Stokes (NS for short) equations. As in the incompressible NS equations, the question of the regularity and uniqueness of weak solutions remains one of the biggest open problems in mathematical fluid mechanics. Energy equality is clearly a prerequisite for regularity, and can be a first step in proving conditional regularity results \cite{SS}.  Prodi \cite{GP},  Serrin \cite{SJ} and Escauriaza-Seregin-Sverak \cite{ELSGA} first proved conditional regularity result
\begin{equation}\label{1.6}
u\in L^{q}\left(0,T;L^{p}\left(\Omega\right)\right) \quad \text{with} \quad\frac{2}{q}+\frac{3}{p}\leq1 \quad \text{and}\quad 3\leq p\leq\infty,
\end{equation}
naturally, energy is conserved under this condition.  However,
Lions \cite{LIONS} proved  that energy equality holds for weak solutions $u\in L^4\left(0, T; L^4(\mathbb{T}^d)\right)$. A few years later, Shinbrot \cite{SHIN} improved upon this result to
\begin{equation}\label{1.7}
u\in L^{q}\left(0,T;L^{p}(\mathbb{T}^d)\right),\,\,  \frac{2}{q}+\frac{2}{p}\leq 1\,\, \text{with}\,\, p\geq 4.
\end{equation}
Note that the energy equality in limiting
case (i.e., $u\in L^{2}\left(0,T;L^{\infty}\right)$) was generalized to $u\in L^{2}\left(0,T; BMO\right)$ by Kozono-Taniuchi \cite{KT}. The recent work by Leslie-Shvydkoy \cite{LS} established the energy
equality under new $L^{q}\left(0,T;L^{p}\right)$ conditions using local energy estimates. In another paper, they \cite{LS1} proved that any solution to the NS equations in $\mathbb{R}^{3}$ which is Type I (i.e. $\left.\|u(t)\|_{L^{\infty}} \leq \frac{C}{\sqrt{T-t}}\right)$ in time must satisfy the energy equality at the first blowup time $T$. Cheskidov-Luo \cite{CLX} showed that the energy equality is holding in the weak-in-time Onsager spaces, as a corollary, they deduced that $u \in L^{q, \infty}\left([0, T] ; B_{p, \infty}^{0}\right)$ for $\frac{1}{q}+\frac{1}{p}=\frac{1}{2},\,p>4$ implies energy equality. It is worth pointing out that the methods developed in \cite{LS1} and \cite{CLX} can not deal with the energy conservation problem in domains with boundaries. For Shinbrot's result, recently, Yu \cite{YUCHENG} gave a new proof which relies on a crucial lemma introduced by Lions. All these results deal with either $ \Omega= \mathbb{R}^d$ or $\Omega= \mathbb{T}^d$. However, due to the well-recognized dominant role of the boundary in the generation of turbulence, it seems very reasonable to investigate the energy conservation in bounded domains. Cheskidov-Friedlander-Shvydkoy \cite{CF} proved energy equality for $u\in L^3D\left(A^\frac{5}{12} \right)$ on a bounded domain, here $A$ denotes the Stokes
operator. Recently, Tan-Yin \cite{TAN}  proved the energy equality holds if $u\in L^{2,\infty}\left(0,T; BMO\right)$, in particular, the methods developed in \cite{TAN} can deal with the local energy conservation problem in domains with boundaries. Berselli and Chiodaroli \cite{BLC} established some new energy balance criteria involving the gradient of the
velocity. Yu \cite{YC} proved the validity of energy equality under the following assumption
$$
u \in L^{q}\left([0, T] ; L^{p}(\Omega)\right) \cap L^{s}\left(0, T ; B_{s}^{\alpha, \infty}(\Omega)\right)
$$
for $\frac{1}{q}+\frac{1}{p} \leq \frac{1}{2}, p \geq 4$ and $\frac{1}{2}+\frac{1}{s}<\alpha<1, s>2$. This result was improved by Nguyen-Nguyen-Tang \cite{NQH} by showing the Shinbrot condition \eqref{1.7} together with a boundary layer assumption $P\in L^2([0,T];L^2(\Omega\setminus\Omega_\delta))$ implies energy conservation. Very recently, Chen-Liang-Wang-Xu \cite{CLWX} showed that the Shinbrot's condition \eqref{1.7} together with a boundary condition $P \in L^{2}\left(0, T ; L^{2}(\partial \Omega)\right)$ guaranteed the energy equality. The additional assumptions $u \in L^{s}\left(0, T ; B_{s}^{\alpha, \infty}(\Omega)\right)$ in \cite{YC}, $P\in L^2([0,T];L^2(\Omega\setminus\Omega_\delta))$ in \cite{NQH} and $P \in L^{2}\left(0, T ; L^{2}(\partial \Omega)\right)$ in \cite{CLWX} are used to deal with the boundary effects, furthermore, they must assume that $\Omega$ is a bounded domain with $C^2$ boundary $\partial \Omega$.


To the best of our knowledge, the results about the energy conservation for the viscoelastic fluids model are very few. Our first aim is to establish the Shinbrot type criteria for energy equality for system \eqref{1.1} in domains with boundaries and answer a problem left in \cite{HZ}.
The second aim concerns the problem of energy conservation for distributional solutions to system \eqref{1.1}, see Definition \ref{def1} for the precise definition. In this case, there is not any available regularity on $(u, \bf F)$, apart from the solution being in
$L^2_{loc}\left((0,T)\times\mathbb{R}^d\right)$. The interest in distributional solutions dates back to
Foias \cite{FC}, who proved their uniqueness under the condition \eqref{1.6}. Later
Fabes, Jones, and Riviere~\cite{FJR} proved the existence of distributional solutions
for the Cauchy problem, while the case of the initial-boundary value problem has been
studied mainly starting from the work of Amann \cite{AH}. As usual, when dealing with
distributional solutions, a duality argument can be employed to show uniqueness, by using
properties of the adjoint problem (which in the case of the Navier-Stokes equations, is a
backward Oseen type-problem). The connection with the non-uniqueness and energy equality has been very
recently studied by Galdi \cite{GALDI,GALDI2} and Buckmaster-Vicol \cite{Vicol}. In \cite{Vicol}, Buckmaster and Vicol showed the non-uniqueness of distributional solutions in $C([0,T], H^\beta(\mathbb{T}^3))$ by constructing energy profile $e(t)=\int_{\mathbb{T}^3}|u(t,x)|^2dx$ satisfying $0=e(t_1)<e(t_2)$ for some $t_1<t_2$. Galdi \cite{GALDI} proved that this non-uniqueness caused by the increase of energy can be inhibited in $L^4L^4$ by showing the distributional solution in $L^4L^4$ is actually a weak solution and the energy equality is valid. It is relevant to observe that the duality argument in \cite{GALDI} is
used to prove that the distributional solution in $L^4L^4$ is the weak limit of a sequence of approximating solutions in $L^\infty L^2\cap L^2H^1$. This also has been used in different
contexts in \cite{BG}, but the approach we follow here takes also inspiration from a
bootstrap argument as used in Lions and Masmoudi \cite{LM} for results concerning mild solutions.

These novelties of our results come from three aspects. Firstly, for the energy equality of the weak solution of \eqref{1.1} on bounded Lipschitz domains, with the help of important properties of weak solutions to the nonstationary
Stokes system and the separate mollification of weak solutions from the boundary effect
to inhibit the boundary effects and prove that although the boundary effect appears, the Shinbrot's criterion on $(u, {\bf F})$ still guarantees the validity of energy conservation of the Leray-Hopf weak solutions, in particular, no boundary layer assumptions are required. This result solves a question left in \cite{HZ}. Secondly, for the energy equality for the weak solution of \eqref{1.1} in whole space $\mathbb{R}^d\, (d=2, 3, 4)$, we show that the behavior of solutions in the finite regions and the behavior at infinite play different roles in energy conservation. This result improves the classical result of Shinbrot \cite{SHIN} by taking $\bf F=0$. Thirdly, for the energy conservation of the distributional solution to system \eqref{1.1}, due to the hyperbolic nature of the equation-{\bf F} of \eqref{1.1},  it seems that the method of Galdi \cite{GALDI} is hard to apply directly to the system \eqref{1.1}.  We provided a simpler and more direct approach which is based on the uniform $L^\infty L^2\cap L^2H^1$ estimates for the approximating solutions $u^\varepsilon=u\ast\eta_\varepsilon$ to establish the energy equality for the distributional solutions to system \eqref{1.1}.

Now, we recall the definitions of Leray-Hopf weak solutions and distributional solutions.

\begin{definition}[Leray-Hopf weak solutions]\label{def1}
Let $(u_0, {\bf F}_{0})\in L^2(\Omega)$ with $\nabla\cdot u_0 = 0$, $T > 0$. The function $(u, {\bf F})$ defined in $[0, T ] \times\Omega$ is said to be a Leray-Hopf weak solution to \eqref{1.1} if\\
1. $u \in L^{\infty}\left(0, T ; L^{2}\left(\Omega\right)\right) \cap L^{2}\left(0, T ; H^{1}\left(\Omega\right)\right), {\bf F} \in L^{\infty}\left(0, T ; L^{2}\left(\Omega\right)\right)$; \\
2. for any $t \in[0, T], \Phi, \Psi \in C_{c}^{\infty}(\Omega \times[0, T])$, we have
$$
\begin{aligned}
	&\int_{\Omega} u(x, t) \cdot \Phi(x, t) d x-\int_{\Omega} u(x, 0) \cdot \Phi(x, 0) d x-\int_{0}^{t} \int_{\Omega} u(x, s) \cdot \partial_{s} \Phi(x, s) d x d s \\
	=&\int_{0}^{t} \int_{\Omega}(u \otimes u: \nabla \Phi- {\bf F} {\bf F^{\top}}: \nabla \Phi+P \dive {\Phi}) d x d s- \int_{0}^{T} \int_{\Omega} \nabla u: \nabla \Phi d x d s
\end{aligned}
$$
and
$$
\begin{aligned}
	&\int_{\Omega} {\bf F}(x, t) \cdot \Psi(x, t) d x-\int_{\Omega} {\bf F}(x, 0) \cdot \Psi(x, 0) d x-\int_{0}^{t} \int_{\Omega} {\bf F}(x, s) \cdot \partial_{s} \Psi(x, s) d x d s \\
	&=\int_{0}^{t} \int_{\Omega}({\bf F} \otimes u: \nabla \Psi-u \otimes {\bf F}: \nabla \Psi) d x d s
\end{aligned}
$$
3. for any $\varphi \in C_{c}^{\infty}(\Omega)$, it holds that
$$
\int_{\Omega} u \cdot \nabla \varphi d x=0
$$
a.e. $t \in(0, T)$;\\
4. $(u,{\bf F})$ satisfies the energy inequality
$$	\frac{1}{2}\int_{\Omega}\left(|u|^{2}+|{\bf F}|^{2}\right)(t) d x+\int_{0}^{t} \int_{\Omega}|\nabla {\bf F}|^{2} d x d s\leq\frac{1}{2} \int_{\Omega}\left(\left|u_{0}\right|^{2}+\left|{\bf F}_{0}\right|^{2}\right) d x,$$
for all $t \in[0, T]$. Furthermore, it holds that $$\lim_{t\to 0^+}\left(\|u(\cdot,t)\|^2_{L^2}+\|{\bf F}(\cdot,t)\|^2_{L^2}\right)=\|u_0\|^2_{L^2}+\|{\bf F}_0\|^2_{L^2}.$$
\end{definition}

In a same fashion with \cite{GALDI}, we define a distributional solution as follows.
\begin{definition}[distributional solutions]
	\label{de2}
	Let $(u_0, {\bf F}_{0})\in L^2(\mathbb{R}^d)$ with $\nabla\cdot u_0 = 0$, $T > 0$. The function $(u, {\bf F})\in L_{\mathrm{loc}}^{2}\left(\mathbb{R}^{d} \times[0, T)\right)$ is a distributional solution to the viscoelastic fluids problem \eqref{1.1} if\\
	1. for any $\Phi \in \mathcal{D}_{T}$ and  $\mathcal{D}_{T}:=\left\{\Phi \in C_{0}^{\infty}\left(\mathbb{R}^{d} \times[0, T)\right): \operatorname{div} \Phi=0\right\}$, we have
	$$
	\begin{aligned}
		&\int_0^T\int_{\mathbb{R}^{d}} u \cdot \partial_{t} \Phi+u \cdot \Delta \Phi+u \otimes u: \nabla \otimes \Phi-{\bf F} {\bf F^{\top}}: \nabla \otimes \Phi d x dt=-\int_{\mathbb{R}^{d}} u(x, 0) \cdot \Phi(x, 0) d x, \\
		&\int_0^T\int_{\mathbb{R}^{d}}  {\bf F} \cdot \partial_{t} \Phi+u \otimes {\bf F}: \nabla \otimes \Phi-{\bf F} \otimes u: \nabla \otimes \Phi d x dt=-\int_{\mathbb{R}^{d}} {\bf F}(x, 0) \cdot \Phi(x, 0) d x;
	\end{aligned}
	$$
	3. for any $\varphi \in C_{0}^{\infty}(\mathbb{R}^{d})$, it holds that
	$$
	\int_{\mathbb{R}^{d}} u \cdot \nabla \varphi d x=0
	$$
	a.e. $t \in(0, T)$;\\
\end{definition}

Now we state our main results.
\begin{theorem}\label{th1}Let $\Omega$ be a bounded domain in $\mathbb{R}^{d}\,(d\geq 2)$ with Lipschitz boundary $\partial\Omega$ and $(u, {\bf F})$ be a weak solution to the initial-boundary value problem \eqref{1.1} in the sense of Definition \ref{def1}. Assume that for any $0<\tau< T$
\begin{equation}\label{1.8}
	u\in L^{q}\left(\tau,T;L^{p}\left(\Omega\right)\right) \quad \text{with} \quad\frac{2}{q}+\frac{2}{p}\leq1, \quad  p\geq 4
\end{equation}
and
\begin{equation}\label{1.9}
	{\bf F}\in L^{r}\left(\tau,T;L^{s}\left(\Omega\right)\right) \quad \text{with} \quad\frac{2}{r}+\frac{2}{s}\leq1, \quad  s\geq 4,
\end{equation}
then the energy equality \eqref{1.2} holds for all $t \in  [0, T ]$.
\end{theorem}
\begin{remark}
The authors in \cite{HZ} proposed that it would be very interesting to investigate energy identity of system \eqref{1.1} with boundary effect. Theorem \ref{th1} gives a positive answer to it. In addition, we don't need to assume $\dive{\bf F^{\top}}= 0$ in the system \eqref{1.1}.
\end{remark}
\begin{remark}
We will note that $$u \in L^{\infty}\left(0, T ; L^{2}(\Omega)\right) \cap L^{q}\left(\tau, T ; L^{p}(\Omega)\right)\text{ for any } \frac{1}{q}+\frac{1}{p} \leq \frac{1}{2},\,\, p \geq 4,$$
and
$$	{\bf F} \in L^{\infty}\left(0, T ; L^{2}(\Omega)\right) \cap L^{r}\left(\tau, T ; L^{s}(\Omega)\right) \text{ for any } \frac{1}{r}+\frac{1}{s} \leq \frac{1}{2},\,\,s \geq 4,$$
then one can deduce that
$$
\|u\|_{L^{4}\left(\tau, T ; L^{4}(\Omega)\right)} \leq C\|u\|_{L^{\infty}\left(\tau, T ; L^{2}(\Omega)\right)}^{a}\|u\|_{L^{q}\left(\tau, T ; L^{p}(\Omega)\right)}^{1-a}
$$
and
$$
\|	{\bf F}\|_{L^{4}\left(\tau, T ; L^{4}(\Omega)\right)} \leq C\|	{\bf F}\|_{L^{\infty}\left(\tau, T ; L^{2}(\Omega)\right)}^{a}\|	{\bf F}\|_{L^{r}\left(\tau, T ; L^{s}(\Omega)\right)}^{1-a}
$$
for some $0<a<1$. Thus, we can use the facts that $u, 	{\bf F}$ are bounded in $L^{4}\left(\tau, T ; L^{4}(\Omega)\right)$ in our proof.
\end{remark}

We can also apply our methods to deal with the case $\Omega= \mathbb{R}^d$, and give an improvement for the result of He and Zi \cite{HZ}.
\begin{theorem}\label{th2} Let $(u, {\bf F}, P)$ be a Leray-Hopf weak solution to the problem \eqref{1.1} in whole space $\mathbb{R}^{d}$ with $d=2, 3, 4$. In
addition, if for any $0<\tau< T$
	\begin{equation}\label{1.11}
		u\in L^{q}\left(\tau,T;L_{loc}^{p}\left(\mathbb{R}^{d}\right)\right) \quad \text{with} \quad\frac{2}{q}+\frac{2}{p}\leq1, \quad  p\geq 4
	\end{equation}
	and
	\begin{equation}\label{1.12}
		{\bf F}\in L^{r}\left(\tau,T;L_{loc}^{s}\left(\mathbb{R}^{d}\right)\right)\cap L^{\frac{4d+8}{d+4}}\left(0,T;L^{\frac{4d+8}{d+4}}\left(\mathbb{R}^{d}\right)\right) \quad \text{with} \quad\frac{2}{r}+\frac{2}{s}\leq1, \quad  s\geq 4,
	\end{equation}
	then the energy equality \eqref{1.2} holds for all $t \in  [0, T ]$.
\end{theorem}
\begin{remark}\label{RE1.3}
From the Theorem \ref{th2}, it seems that the behavior of the solution on the finite regions and the behavior at infinite play different roles in the study of energy conservation to problem \eqref{1.1}.
\end{remark}
Next, the objective of the following theorem is to prove that, actually, for the validity of energy equality, the requirement  $u \in L^{\infty}\left(0, T ; L^{2}\left(\mathbb{R}^{d}\right)\right) \cap L^{2}\left(0, T ; H^{1}\left(\mathbb{R}^{d}\right)\right)$ and $ {\bf F} \in L^{\infty}\left(0, T ; L^{2}\left(\mathbb{R}^{d}\right)\right)$ can be removed. It is just enough that $(u,{\bf F}) \in L^{4}\left(0, T ; L^{4}(\mathbb{R}^{d})\right)$, along with the (necessary) condition $(u_0, {\bf F}_0)\in L^2(\mathbb{R}^d)$.
\begin{theorem}\label{th3}
	Let $(u, {\bf F})\in L_{\mathrm{loc}}^{2}\left(\mathbb{R}^{d} \times[0, T)\right)\,(d\geq 2)$ be a distributional solution in the sense of Definition \ref{de2} to the system \eqref{1.1}. In additional, if $(u,{\bf F}) \in L^{4}\left(0, T ; L^{4}(\mathbb{R}^{d})\right)$, then
	\begin{equation*}
		\begin{split}
			\int_{\mathbb{R}^d}|u(x,t)|^{2}+|{\bf F}(x,t)|^{2} d x+2 \int_{0}^{t} \int_{\mathbb{R}^{d}}|\nabla u|^{2} d x d t=\int_{\mathbb{R}^{d}}\left|u_{0}\right|^{2}+\left|{\bf F}_{0}\right|^{2} d x
		\end{split}
	\end{equation*}
	for any $t \in[0, T]$.
\end{theorem}
\begin{remark}\label{re1}
	This result extends the well-known Galdi's energy conservation criterion $u\in L^4L^4$ for a distributional solution to the incompressible viscoelastic fluids \eqref{1.1}.
\end{remark}
\section{Preliminaries}
To investigate the boundary effect, we will recall some useful notations and Lemmas.\\
Assume $\Omega$ be a bounded domain with boundary $\partial \Omega$ in $\mathbb{R}^{d}$ with $d \geq 2$, for any $x\in\Omega$, $d(x)=\inf_{y\in\partial\Omega}|x-y|$, we define $\Omega^*_{\delta}$ as follows
$$
\Omega^*_{\delta}=\{x\in\Omega|d(x)>\delta\}.
$$
It is clear that $\Omega^*_{\delta}$ is still a domain if $\delta$ is small enough. We use $L^{q}\left(0, T ; L^{p}(\Omega)\right)$ to denote the space of measurable functions with the following norm
$$
\|f\|_{L^{q}\left(0, T ; L^{p}(\Omega)\right)}=\left\{\begin{array}{l}
	\left(\int_{0}^{T}\left(\int_{\Omega}|f(t, x)|^{p} d x\right)^{\frac{q}{p}} d t\right)^{\frac{1}{q}}, 1 \leq q<\infty \\
	\operatorname{ess} \sup _{t \in[0, T]}\|f(t, \cdot)\|_{L^{p}(\Omega)}, q=\infty.
\end{array}\right.
$$
 If $u\in L^q(\tau,T;L^p(\Omega))$ for any $0<\tau<T$, we say $u\in L^q_{loc}(0,T;L^p(\Omega))$.

Let $\eta: \mathbb{R}^{d} \rightarrow \mathbb{R}$ be a standard mollifier, i.e. $\eta(x)=C \mathrm{e}^{\frac{1}{|x|^{2}-1}}$ for $|x|<1$ and $\eta(x)=0$ for $|x| \geqslant 1$, where  constant $C>0$ selected such that $\int_{\mathbb{R}^{d}} \eta(x) \mathrm{d} x=1$. For any $\varepsilon>0$, we define the rescaled mollifier $\eta_{\varepsilon}(x)=\varepsilon^{-d} \eta\left(\frac{x}{\varepsilon}\right) .$ For any function $f \in L_{\text {loc }}^{1}(\Omega)$, its mollified version is defined as
$$f^{\varepsilon}(x)=\left(f * \eta_{\varepsilon}\right)(x)=\int_{\Omega} \eta_{\varepsilon}(x-y) f(y) \mathrm{d} y .$$
If $f \in  W^{1, p}(\Omega)$, the following local approximation is well known
$$
f^{\varepsilon}(x) \rightarrow f \quad \text { in } \quad  W_{l o c}^{1, p}(\Omega) \quad \forall p \in[1, \infty).
$$
The crucial ingredients to prove main results are the following important lemmas. In the first lemma, we shall construct some cut-off functions to separate the mollification of the weak solution from the boundary effect.
\begin{lemma}\label{CUT-OFF}
Let $\Omega$ be a bounded domain with boundary $\partial\Omega$. For any $\delta>0$ small enough, there exist some cut-off functions $\varphi(x)$ which satisfy $\varphi(x)=1$ for $x \in \Omega^*_{3 \delta}$ and $\varphi(x)=0$ for $x \in \Omega \backslash \Omega^*_{\delta}$. Furthermore, it holds that
$$
|\nabla \varphi(x)| \leq \frac{C}{\delta}.
$$
\textbf{Proof.} We first take $\chi(x)$ is the characteristic function of $\Omega^*_{2 \delta}$. We claim that for any $x \in \Omega^*_{2 \delta}$, it yields that
$$
\inf _{y \in \partial \Omega^*_{\delta}}|x-y| \geq \delta.
$$
We thus define $\varphi(x)= (\chi*\eta_{\delta})(x)$, it is obvious that $\varphi(x)=1$ for $x \in \Omega^*_{3 \delta}$ and $\varphi(x)=0$ for $x \in \Omega \backslash \Omega^*_{\delta}$
and $|\nabla \varphi(x)| \leq \frac{C}{\delta}$. We now prove our claim, for any $x \in \Omega^*_{2 \delta}$, assume $y_{0} \in \partial \Omega^*_{\delta}$ satisfy $\left|x-y_{0}\right|=$ $\inf _{y \in \partial \Omega^*_{\delta}}|x-y| .$ Similarly, we assume $z_{0} \in \partial \Omega$ satisfy $\left|y_{0}-z_{0}\right|=\inf _{z \in \partial \Omega}\left|y_{0}-z\right| .$ By the definitions
of $\Omega^*_{\delta}$ and $\Omega^*_{2 \delta}$, we have
$$
\left|x-z_{0}\right| \geq 2 \delta,\left|y_{0}-z_{0}\right|=\delta.
$$
This means that $\inf _{y \in \partial \Omega^*_{\delta}}|x-y|=\left|x-y_{0}\right| \geq\left|x-z_{0}\right|-\left|y_{0}-z_{0}\right| \geq \delta$.
\end{lemma}
\begin{lemma}\label{le2.2} \cite{NQH}
Let $2 \leq d \in \mathbb{N}, \Omega\subset\mathbb{R}^d$ be a bounded domain with boundary $\partial\Omega$, $1 \leq p, q \leq \infty$ and $f, g: \Omega \times(0, T) \rightarrow \mathbb{R}$.

	(i) Assume $f \in L^{p}\left(0, T ; L^{q}\left(\Omega\right)\right)$. Then for any $0<\varepsilon<\delta$, there holds
	$$
	\left\|\nabla f^{\varepsilon}\right\|_{L^{p}\left(0, T; L^{q}\left(\Omega^*_\delta\right)\right)} \leq C \varepsilon^{-1}\|f\|_{L^{p}\left(0, T; L^{q}\left(\Omega\right)\right)}
	$$
	Moreover, if $p, q<\infty$ then
	$$
	\underset{\varepsilon \rightarrow 0}\limsup \varepsilon\left\|\nabla f^{\varepsilon}\right\|_{L^{p}\left(0, T; L^{q}\left(\Omega^*_\delta\right)\right)}=0.
	$$
	
	(ii) Let $p, p_{1} \in[1, \infty)$ and $p_{2} \in(1, \infty]$ with $\frac{1}{p}=\frac{1}{p_{1}}+\frac{1}{p_{2}} .$ Assume $f \in L^{p_{1}}\left(0, T ; W^{1, p_{1}}\left(\Omega\right)\right)$ and $g \in L^{p_{2}}\left(\Omega \times(0, T)\right)$. Then for any $\frac{\delta}{2}>\varepsilon>0$, there holds
	$$
	\left\|(f g)^{\varepsilon}-f^{\varepsilon} g^{\varepsilon}\right\|_{L^{p}\left(\Omega^*_{\delta} \times(0, T)\right)} \leqslant C \varepsilon\|f\|_{L^{p_{1}}\left(0, T ; W^{1, p_{1}}\left(\Omega^*_{\frac{\delta}{2}}\right)\right)}\|g\|_{L^{p_{2}}\left(\Omega^*_{\frac{\delta}{2}} \times(0, T)\right)}
	$$
	Moreover, if $p_{2}<\infty$ then
	$$
	\underset{\varepsilon \rightarrow 0}{\limsup } \varepsilon^{-1}\left\|(f g)^{\varepsilon}-f^{\varepsilon} g^{\varepsilon}\right\|_{L^{p}\left(\Omega^*_{\delta}  \times(0, T)\right)}=0.
	$$
	\end{lemma}
Similar results can be obtained for $\Omega=\mathbb{R}^d$ using only minor modifications.
\begin{lemma}\label{le2.3} \cite{NQH}
	Let $2 \leq d \in \mathbb{N}$, $1 \leq p, q \leq \infty$ and $f, g: \mathbb{R}^d \times(0, T) \rightarrow \mathbb{R}$.
	
	(i) Assume $f \in L^{p}\left(0, T ; L^{q}\left(\mathbb{R}^d\right)\right)$. Then for any $0<\varepsilon<\delta$, there holds
	$$
	\left\|\nabla f^{\varepsilon}\right\|_{L^{p}\left(0, T; L^{q}\left(\mathbb{R}^d\right)\right)} \leq C \varepsilon^{-1}\|f\|_{L^{p}\left(0, T; L^{q}\left(\mathbb{R}^d\right)\right)}
	$$
	Moreover, if $p, q<\infty$ then
	$$
	\underset{\varepsilon \rightarrow 0}\limsup \varepsilon\left\|\nabla f^{\varepsilon}\right\|_{L^{p}\left(0, T; L^{q}\left(\mathbb{R}^d\right)\right)}=0.
	$$
	
	(ii) Let $p, p_{1} \in[1, \infty)$ and $p_{2} \in(1, \infty]$ with $\frac{1}{p}=\frac{1}{p_{1}}+\frac{1}{p_{2}} .$ Assume $f \in L^{p_{1}}\left(0, T ; W^{1, p_{1}}\left(\mathbb{R}^d\right)\right)$ and $g \in L^{p_{2}}\left(\mathbb{R}^d \times(0, T)\right)$. Then for any $\delta>\varepsilon>0$, there holds
	$$
	\left\|(f g)^{\varepsilon}-f^{\varepsilon} g^{\varepsilon}\right\|_{L^{p}\left(\mathbb{R}^d \times(0, T)\right)} \leqslant C \varepsilon\|f\|_{L^{p_{1}}\left(0, T ; W^{1, p_{1}}\left(\mathbb{R}^d\right)\right)}\|g\|_{L^{p_{2}}\left(\mathbb{R}^d \times(0, T)\right)}
	$$
	Moreover, if $p_{2}<\infty$ then
	$$
	\underset{\varepsilon \rightarrow 0}{\limsup }   \varepsilon^{-1}\left\|(f g)^{\varepsilon}-f^{\varepsilon} g^{\varepsilon}\right\|_{L^{p}\left(\mathbb{R}^d  \times(0, T)\right)}=0.
	$$
\end{lemma}
\begin{lemma}\label{le2.4}\cite{SOHR}
	Let $\Omega \subseteq \mathbb{R}^{d}, d \geq 2$, be any domain, let $0<T \leq \infty$, $u_{0} \in$ $L_{\sigma}^{2}(\Omega), f=f_{0}+\operatorname{div} F$ with
	$$
	f_{0} \in L^{1}\left(0, T ; L^{2}(\Omega)\right), \quad F \in L^{2}\left(0, T ; L^{2}(\Omega)\right),
	$$
	and let
	$$
	u \in L_{l o c}^{1}\left([0, T) ; W_{0, \sigma}^{1,2}(\Omega)\right)
	$$
	be a weak solution of the Stokes system
	$$
	u_{t}-\mu \Delta u+\nabla P=f, \quad \operatorname{div} u=0,\left.\quad u\right|_{\partial \Omega}=0, \quad u(0)=u_{0}
	$$
	with data $f, u_{0}$.
	Then $u$ has the following properties:
	
	a) $u \in L^{\infty}\left(0, T ; L_{\sigma}^{2}(\Omega)\right), \quad \nabla u \in L^{2}\left(0, T ; L^{2}(\Omega)\right) .$
	
	b) $u:[0, T) \rightarrow L_{\sigma}^{2}(\Omega)$ is strongly continuous, after a redefinition on a null set of $[0, T), u(0)=u_{0}$, and the energy equality
	$$
	\begin{aligned}
		\frac{1}{2}\|u(t)\|_{2}^{2}+\mu \int_{0}^{t}\|\nabla u\|_{2}^{2} d \tau=& \frac{1}{2}\left\|u_{0}\right\|_{2}^{2}+\int_{0}^{t}\langle f_{0}, u\rangle_{\Omega} d \tau -\int_{0}^{t}\langle F, \nabla u\rangle_{\Omega} d \tau.
	\end{aligned}
	$$	
\end{lemma}
\section{On the energy equality for Leray-Hopf weak solutions.}
\subsection{Proof of Theorem \ref{th1}}
The main object of this subsection is to prove Theorem \ref{th1}. Let us first rewrite the first equation of the system \eqref{1.1} as
$$
\partial_t u-\Delta u+\nabla P=\nabla\cdot\left(\left({\bf F}{\bf F}^\top\right)-\left(u\otimes u\right)\right),
$$
due to
$$\left(\left({\bf F}{\bf F}^\top\right)-(u\otimes u)\right) \in L_{l o c}^{2}\left((0, T) ; L^{2}(\Omega)\right),$$
it follows from Lemma \ref{le2.4} that $u$ sastisfies
\begin{equation}
	\begin{aligned}\label{ST}
		\frac{1}{2}\|u(t)\|_{L^2}^{2}+\int_{0}^{t}\|\nabla u\|_{L^2}^{2} d \tau= \frac{1}{2}\left\|u_{0}\right\|_{L^2}^{2} -\int_{0}^{t}\langle ({\bf F}{\bf F}^\top)-(u\otimes u), \nabla u\rangle_{\Omega} d \tau.
	\end{aligned}
\end{equation}
Next, for the sake of simplicity, we will proceed as if the
solution is differentiable in time. The extra arguments needed to mollify in time are straightforward. We mollify the second equation of system \eqref{1.1},  then using  ${\bf F}^{\varepsilon}(t, x) \varphi(x)$ to test the resulting equation, integrating by parts over $[\tau, t] \times \Omega$ with $0<\tau \leq t \leq T$, one has
\begin{equation}\label{A}
\begin{aligned}
	&\frac{1}{2} \int_{\Omega}\left|{\bf F}^{\varepsilon}(t, x)\right|^{2} \varphi(x) d x-\frac{1}{2} \int_{\Omega}\left|{\bf F}^{\varepsilon}(\tau, x)\right|^{2} \varphi(x) d x\\
	=&\int_{\tau}^{t} \int_{\Omega}(\nabla u{\bf F})^{\varepsilon}{\bf F}^{\varepsilon}\varphi(x) d x d s-\int_{\tau}^{t} \int_{\Omega}(u\cdot\nabla {\bf F})^{\varepsilon}{\bf F}^{\varepsilon}\varphi(x) d x d s,
\end{aligned}
\end{equation}
where $\varphi(x)$ is a cut-off function constructed in Lemma \ref{CUT-OFF} and which equals to one on $\Omega^*_{3 \delta}$ and vanishes out of $\Omega^*_{\delta}$. We claim that
\begin{equation}\label{B}
	\int_{\tau}^{t} \int_{\Omega}(\nabla u{\bf F})^{\varepsilon}{\bf F}^{\varepsilon}\varphi(x) d x d s\rightarrow \int_{\tau}^{t} \int_{\Omega}(\nabla u{\bf F}){\bf F}\varphi(x) d x d s,\,\, as \,\,\varepsilon \rightarrow 0.
\end{equation}
Indeed,
\begin{equation*}
	\begin{split}
		&\left|\int^t_\tau\int_{\Omega}(\nabla u{\bf F})^{\varepsilon} {\bf F}^{\varepsilon}\varphi-(\nabla u{\bf F}) {\bf F}\varphi  dxds\right|\\
		=&\left|\int^t_\tau\int_{\Omega}(\partial_ku_i{\bf F}_{kj})^{\varepsilon} {\bf F}_{ij}^{\varepsilon} \varphi-(\partial_ku_i{\bf F}_{kj}) {\bf F}_{ij} \varphi  dxds\right|\\
		=&\left|\int^t_\tau\int_{\Omega}[(\partial_ku_i{\bf F}_{kj})^{\varepsilon} -(\partial_ku_i{\bf F}_{kj})]{\bf F}_{ij}^{\varepsilon} \varphi+(\partial_ku_i{\bf F}_{kj}){\bf F}_{ij}^{\varepsilon} \varphi-(\partial_ku_i{\bf F}_{kj}) {\bf F}_{ij} \varphi  dxds\right|\\
		\leq&C\|(\partial_ku_i{\bf F}_{kj})^{\varepsilon} -(\partial_ku_i{\bf F}_{kj})\|_{L^{\frac{4}{3}}\left(\tau, t ; L^{\frac{4}{3}}(\Omega)\right)} \|{\bf F}_{ij} \varphi \|_{L^{4}\left(\tau, t ; L^{4}(\Omega)\right)}\\
		&+C\|{\bf F}_{ij}^{\varepsilon} \varphi-{\bf F}_{ij} \varphi\|_{L^{4}\left(\tau, t ; L^{4}(\Omega)\right)} \|\partial_ku_i{\bf F}_{kj}\|_{L^{\frac{4}{3}}\left(\tau, t ; L^{\frac{4}{3}}(\Omega)\right)}\\
		&\rightarrow 0,\,\, as \,\,\varepsilon \rightarrow 0.
	\end{split}
\end{equation*}
We next rewrite
\begin{equation}\label{C}
	\begin{aligned}
		&\int_{\tau}^{t} \int_{\Omega}(u\cdot\nabla {\bf F})^{\varepsilon} {\bf F}^{\varepsilon}\varphi d x d s \\
		=&\int_{\tau}^{t} \int_{\Omega}(u_k\partial_k {\bf F}_{ij})^{\varepsilon} {{\bf F}^{\varepsilon}_{ij}}\varphi d x d s\\
		=&\int_{\tau}^{t} \int_{\Omega}\partial_k(u_k {\bf F}_{ij})^{\varepsilon} {{\bf F}^{\varepsilon}_{ij}}\varphi d x d s\\
		=&-\int_{\tau}^{t} \int_{\Omega}(u_k {\bf F}_{ij})^{\varepsilon} \partial_k({{\bf F}^{\varepsilon}_{ij}}\varphi)  d x d s\\
		=&-\int_{\tau}^{t} \int_{\Omega}[(u_k {\bf F}_{ij})^{\varepsilon}-(u^{\varepsilon}_k {\bf F}^{\varepsilon}_{ij})] \partial_k({{\bf F}^{\varepsilon}_{ij}}\varphi) dxds-\frac{1}{2}\int_{\tau}^{t} \int_{\Omega}u^{\varepsilon}_k\partial_k\varphi|{{\bf F}^{\varepsilon}_{ij}}|^2
		dxds\\
		\leq &C\|(u_k {\bf F}_{ij})^{\varepsilon}-(u^{\varepsilon}_k {\bf F}^{\varepsilon}_{ij})\|_{L^{\frac{4}{3}}\left(\tau, t ; L^{\frac{4}{3}}(\Omega^*_\delta)\right)} \|\partial_k({{\bf F}^{\varepsilon}_{ij}}\varphi)\|_{L^{4}\left(\tau, t ; L^{4}(\Omega^*_\delta)\right)}-\frac{1}{2}\int_{\tau}^{t} \int_{\Omega}u^{\varepsilon}_k\partial_k\varphi|{{\bf F}^{\varepsilon}_{ij}}|^2
		dxds\\
		\leq &C\varepsilon^{-1}\|(u_k {\bf F}_{ij})^{\varepsilon}-(u^{\varepsilon}_k {\bf F}^{\varepsilon}_{ij})\|_{L^{\frac{4}{3}}\left(\tau, t ; L^{\frac{4}{3}}(\Omega^*_\delta)\right)} \|{{\bf F}_{ij}}\|_{L^{4}\left(\tau, t ; L^{4}(\Omega)\right)}-\frac{1}{2}\int_{\tau}^{t} \int_{\Omega}u^{\varepsilon}_k\partial_k\varphi|{{\bf F}^{\varepsilon}_{ij}}|^2
		dxds\\
		&\rightarrow -\frac{1}{2}\int_{\tau}^{t} \int_{\Omega}u_k\partial_k\varphi|{{\bf F}_{ij}}|^2,\,\, as \,\,\varepsilon \rightarrow 0,
	\end{aligned}
\end{equation}
where we used the key Lemma \ref{le2.2}. Hence
\begin{equation}\label{D}
\int_{\tau}^{t} \int_{\Omega}(u\cdot\nabla {\bf F})^{\varepsilon} {\bf F}^{\varepsilon}\varphi d x d s\rightarrow -\frac{1}{2}\int_{\tau}^{t} \int_{\Omega}u_k\partial_k\varphi|{{\bf F}_{ij}}|^2
dxds,\,\, as \,\,\varepsilon \rightarrow 0.
\end{equation}
We first take $\varepsilon \rightarrow 0$ in \eqref{A} and substituting all above results into \eqref{A} gives that
\begin{equation}\label{E}
	\begin{aligned}
		& \frac{1}{2}\int_{\Omega}\left|{\bf F}(t, x)\right|^{2}\varphi(x) d x-\frac{1}{2}\int_{\Omega}\left|{\bf F}(\tau, x)\right|^{2}\varphi(x)dx\\
		=&\frac{1}{2}\int_{
			\tau}^{t} \int_{\Omega}u_k\partial_k\varphi |{\bf F}_{ij}|^2dxds+\int_{\tau}^{t}\int_{\Omega}(\partial_ku_i{\bf F}_{kj}){\bf F}_{ij}\varphi dxds.
	\end{aligned}
\end{equation}
Using Lemma \ref{CUT-OFF} and H\"older's inequality, we have
$$
\begin{aligned}
\frac{1}{2}\int_{
	\tau}^{t} \int_{\Omega}u_k\partial_k\varphi |{\bf F}_{ij}|^2dxds\leq	 \frac{C}{\delta}\|{\bf F}\|_{L^{4}\left(\tau, T ; L^{4}\left(\Omega \backslash \Omega^*_{3 \delta}\right)\right)}^{2}
 \cdot\|u\|_{L^{2}\left(\tau, T ; L^{2}\left(\Omega \backslash \Omega^*_{3 \delta}\right)\right)}.
\end{aligned}
$$
Noticing that $u|_{\partial \Omega}=0$ and the $\Omega$ is a Lipschitz domain, by using Poincar\'e's inequality, we find that
$$
\|u\|_{L^{2}\left(\tau, T ; L^{2}\left(\Omega \backslash \Omega^*_{3 \delta}\right)\right)} \leq C \delta\|\nabla u\|_{L^{2}\left(\tau, T ; L^{2}\left(\Omega \backslash \Omega^*_{3 \delta}\right)\right)}.
$$
Obviously,
$$
\begin{aligned}
	\frac{1}{2}\int_{
		\tau}^{t} \int_{\Omega}u_k\partial_k\varphi |{\bf F}_{ij}|^2dxds
	&\leq	 C\|{\bf F}\|_{L^{4}\left(\tau, T ; L^{4}\left(\Omega \backslash \Omega^*_{3 \delta}\right)\right)}^{2}\cdot\|\nabla u\|_{L^{2}\left(\tau, T ; L^{2}\left(\Omega \backslash \Omega^*_{3 \delta}\right)\right)}\\
	& \rightarrow 0, \quad \text { as } \delta \rightarrow 0.
\end{aligned}
$$
Letting $\delta$ go to zero in \eqref{E}, and using above estimates involving $\frac{1}{2}\int_{
	\tau}^{t} \int_{\Omega}u_k\partial_k\varphi |{\bf F}_{ij}|^2dxds \rightarrow 0$, we obtain
\begin{equation}\label{F}
\frac{1}{2}\int_{\Omega}\left|{\bf F}(t, x)\right|^{2} d x-\frac{1}{2}\int_{\Omega}\left|{\bf F}(\tau, x)\right|^{2} dx-\int_{\tau}^{t}\int_{\Omega}(\partial_ku_i{\bf F}_{kj}){\bf F}_{ij} dxds=0
\end{equation}
for all $0<\tau\leq t\leq T$.

Combining \eqref{ST} and \eqref{F} gives
\begin{equation}
	\begin{aligned}\label{G}
		&\frac{1}{2}\left(\|u(\cdot,t)\|^2_{L^2}
		+\|{\bf F}(\cdot,t)\|^2_{L^2}\right)+\int_{0}^{t}\|\nabla u\|_{L^2}^{2} d s\\
		=& \frac{1}{2}\left(\left\|u_{0}\right\|_{L^2}^{2}+\left\|{\bf F}(\tau)\right\|_{L^2}^{2} \right)-\int_{0}^{t}\langle ({\bf F}{\bf F}^\top)-(u\otimes u), \nabla u\rangle_{\Omega} d s+\int_{\tau}^{t}\int_{\Omega}(\partial_ku_i{\bf F}_{kj}){\bf F}_{ij} dxds.
	\end{aligned}
\end{equation}
By passing to the limit as $\tau\rightarrow 0$ and using the facts that $$\lim_{t\to 0^+}\left(\|u(\cdot,t)\|^2_{L^2}+\|{\bf F}(\cdot,t)\|^2_{L^2}\right)=\|u_0\|^2_{L^2}+\|{\bf F}_0\|^2_{L^2},$$ we immediately get
\begin{equation}
	\begin{aligned}\label{H}
		&\frac{1}{2}\left(\|u(\cdot,t)\|^2_{L^2}
		+\|{\bf F}(\cdot,t)\|^2_{L^2}\right)+\int_{0}^{t}\|\nabla u\|_{L^2}^{2} d s\\
		=& \frac{1}{2}\left(\left\|u_{0}\right\|_{L^2}^{2}+\left\|{\bf F}_0\right\|_{L^2}^{2} \right)+\int_{0}^{t}\langle (u\otimes u), \nabla u\rangle_{\Omega} d s.
	\end{aligned}
\end{equation}
Next, we will prove that
$$\int_{0}^{t} \int_{\Omega}(u\otimes u) \nabla u dx ds\equiv0.$$
Now, $u(t) \in H^{1}_{0}(\Omega)$, for a.a. $t \in[0, T)$ and so, for any such fixed $t$, denoting by $\left\{\psi_{k}\right\}$ a sequence from $\left\{\psi_{k}\right\} \subset\mathcal{D}(\Omega)=\left\{\boldsymbol{\psi} \in C_{0}^{\infty}(\Omega): \operatorname{div} \psi=0 \text { in } \Omega\right\}$ converging to $u$ in $H^{1}$ we have
$$
\begin{aligned}
	\left|\langle u \cdot \nabla u, u\rangle-\langle u \cdot \nabla \psi_{k}, \psi_{k}\rangle\right| & \leq\left|\langle u \cdot \nabla u,\left(u-\psi_{k}\right)\rangle\right|+\left|\langle u \cdot \nabla\left(u-\psi_{k}\right), u\rangle\right| \\
	& \leq\|u\|_{L^4}\|\nabla u\|_{L^2}\left\|u-\psi_{k}\right\|_{L^4}+\|u\|_{L^4}^{2}\left\|\nabla\left(u-\psi_{k}\right)\right\|_{L^2}
\end{aligned}
$$
so we deduce
$$
\lim _{k \rightarrow \infty}\langle u \cdot \nabla \psi_{k}, \psi_{k}\rangle=\langle u \cdot \nabla u, u\rangle.
$$
However, since $\nabla\cdot u(t)=0$ for a.a. $t$, we get
$$
\langle u \cdot \nabla u, u\rangle=\lim_{k\to\infty}\langle u \cdot \nabla \psi_{k}, \psi_{k}\rangle=\lim_{k\to\infty}\frac{1}{2}\langle u, \nabla\left(\psi_{k}\right)^{2}\rangle=0.
$$
Furthermore, $\langle u\cdot \nabla u, u\rangle\in L^1([0,T])$, hence
$$\int_{0}^{t} \int_{\Omega}(u\otimes u) \nabla u dx ds=0.$$
Thus, we immediately get
$$
	\frac{1}{2}\int_{\Omega}\left(|u|^{2}+|{\bf F}|^{2}\right)(t) d x+\int_{0}^{t} \int_{\Omega}|\nabla u|^{2} d x d s=\frac{1}{2} \int_{\Omega}\left(\left|u_{0}\right|^{2}+\left|{\bf F}_{0}\right|^{2}\right) d x.
$$
This completes the proof of Theorem \ref{th1}.
\subsection{Proof of Theorem \ref{th2}}
In this subsection, we shall give the proof of Theorem \ref{th2}. Before going to do it, we first let $\phi \in C_{0}^{\infty}\left(B_{2r}\right)$ be a cut-off function such that $\phi=1$ in $B_{r}$, and
$|\nabla \phi| \leq \frac{C}{r}$ for $r> 0$.
And then,  We mollify the first and the second equation of system \eqref{1.1},  by using  $u^{\varepsilon}(t, x) \phi(x)$ and ${\bf F}^{\varepsilon}(t, x) \phi(x)$ to test the regularized equations of \eqref{1.1}, respectively, after integration by parts over $[\tau, t] \times \mathbb{R}^d$ with $0<\tau \leq t \leq T$, that
\begin{equation}\label{3.19}
	\begin{aligned}
		& \frac{1}{2}\int_{\mathbb{R}^d}\left(\left|u^{\varepsilon}(t, x)\right|^{2}\phi(x)+\left|{\bf F}^{\varepsilon}(t, x)\right|^{2}\phi(x)\right) d x-\frac{1}{2}\int_{\mathbb{R}^d}\left(\left|u^{\varepsilon}(\tau, x)\right|^{2}\phi(x)+\left|{\bf F}^{\varepsilon}(\tau, x)\right|^{2}\phi(x)\right) dx\\
		&+ \int_{\tau}^{t}  \int_{\mathbb{R}^d}\left|\nabla u^{\varepsilon}(s,x)\right|^{2} \phi(x)d x ds\\
		=& \int_{\tau}^{t} \int_{\mathbb{R}^d}\left\{(u \otimes u)^{\varepsilon}: \nabla\left(u^{\varepsilon} \varphi\right)+P^{\varepsilon} \nabla \cdot\left(u^{\varepsilon} \phi\right)-\nabla u^{\varepsilon}:\left(u^{\varepsilon} \otimes \nabla \phi\right)\right\}dxds\\
		&-\int_{\tau}^{t} \int_{\mathbb{R}^d}\left\{({\bf F}{\bf F^{\top}})^{\varepsilon}: \nabla\left(u^{\varepsilon} \phi\right)-(\nabla u{\bf F})^{\varepsilon} {\bf F}^{\varepsilon}\phi+(u\cdot\nabla {\bf F})^{\varepsilon} {\bf F}^{\varepsilon}\phi\right\} d x d s \\
		=& \int_{\tau}^{t} \int_{\mathbb{R}^d}\left\{[(u \otimes u)^{\varepsilon}-(u^{\varepsilon} \otimes u^{\varepsilon})]: \nabla\left(u^{\varepsilon} \phi\right)+(u^{\varepsilon} \otimes u^{\varepsilon}): \nabla\left(u^{\varepsilon} \phi\right)+P^{\varepsilon} \nabla \cdot\left(u^{\varepsilon} \phi\right)-\nabla u^{\varepsilon}:\left(u^{\varepsilon} \otimes \nabla \phi\right)\right\}dxds\\
		&-\int_{\tau}^{t} \int_{\mathbb{R}^d}\left\{({\bf F}{\bf F^{\top}})^{\varepsilon}: \nabla\left(u^{\varepsilon} \phi\right)-(\nabla u{\bf F})^{\varepsilon} {\bf F}^{\varepsilon}\phi+(u\cdot\nabla {\bf F})^{\varepsilon} {\bf F}^{\varepsilon}\phi\right\} d x d s \\
		=& \int_{\tau}^{t} \int_{\mathbb{R}^d}\left\{[(u \otimes u)^{\varepsilon}-(u^{\varepsilon} \otimes u^{\varepsilon})]: \nabla\left(u^{\varepsilon} \phi\right)+\frac{\left|u^{\varepsilon}\right|^{2}}{2} u^{\varepsilon} \cdot \nabla \varphi+P^{\varepsilon} \nabla \cdot\left(u^{\varepsilon} \phi\right)-\nabla u^{\varepsilon}:\left(u^{\varepsilon} \otimes \nabla \phi\right)\right\}dxds\\
		&-\int_{\tau}^{t} \int_{\mathbb{R}^d}\left\{({\bf F}{\bf F^{\top}})^{\varepsilon}: \nabla\left(u^{\varepsilon} \varphi\right)-(\nabla u{\bf F})^{\varepsilon} {\bf F}^{\varepsilon}\varphi+(u\cdot\nabla {\bf F})^{\varepsilon} {\bf F}^{\varepsilon}\phi\right\} d x d s\\
		=: & \mathcal{F}^\varepsilon(t).
	\end{aligned}
\end{equation}
We note that by the construction of the cut-off function $\phi$, this implies that
$$\left|\int_{\tau}^{t} \int_{\mathbb{R}^d}[(u \otimes u)^{\varepsilon}-(u^{\varepsilon} \otimes u^{\varepsilon})]: \nabla\left(u^{\varepsilon} \phi\right) dxds\right |=\left|\int_{\tau}^{t} \int_{B_{2r}}[(u \otimes u)^{\varepsilon}-(u^{\varepsilon} \otimes u^{\varepsilon})]: \nabla\left(u^{\varepsilon} \phi\right) dxds\right |.$$
Due to $u\in L^4\left(0,T;L_{loc}^4(\mathbb{R}^d)\right)$, we get
$$
\begin{aligned}
	&\left|\int_{\tau}^{t} \int_{\mathbb{R}^3}\left[(u \otimes u)^{\varepsilon}-\left(u^{\varepsilon} \otimes u^{\varepsilon}\right)\right]: \nabla\left(u^{\varepsilon} \phi\right) d x d s\right| \\
	\leq &\left(\left\|(u \otimes u)^{\varepsilon}-(u \otimes u)\right\|_{L^{2}\left(\tau, T ; L^{2}\left(B_{2r}\right)\right)}+\left\|\left(u^{\varepsilon} \otimes u^{\varepsilon}\right)-(u \otimes u)\right\|_{L^{2}\left(\tau, T ; L^{2}\left(B_{2r}\right)\right)}\right)\left\|\nabla\left(u^{\varepsilon} \varphi\right)\right\|_{L^{2}\left(\tau, T ; L^{2}\left(B_{2r}\right)\right)}\\
	\rightarrow & 0, \quad as\quad \varepsilon\rightarrow 0.
\end{aligned}
$$
Now, in a similar manner to the convergence results of the mollified terms in Theorem \ref{th1}, when $\varepsilon\rightarrow 0$, we immediately obtain
 $\mathcal{F}^\varepsilon(t)$ have the following limit:
\begin{equation}\label{3.20}
	\begin{aligned}
| \mathcal{F}^\varepsilon(t)|\rightarrow& \bigg|\frac{1}{2}\int_{\tau}^{t} \int_{\mathbb{R}^d}|u|^2u\cdot\nabla\phi dxds+ \int_{\tau}^{t} \int_{\mathbb{R}^d}\left\{P  \left(u\cdot \nabla\phi\right)-\nabla u:\left(u \otimes \nabla \phi\right)\right\}dxds\\
&+\frac{1}{2}\int_{
	\tau}^{t} \int_{\mathbb{R}^d}u_k\partial_k\phi |{\bf F}_{ij}|^2dxds-\int_{\tau}^{t}\int_{\mathbb{R}^d}u_i{\bf F}_{ij}{\bf F}_{kj}\partial_k\phi dxds\bigg|.\\
\end{aligned}
\end{equation}
Letting $\varepsilon \rightarrow 0 $ in \eqref{3.19} and combining \eqref{3.20}, one has
\begin{equation}\label{3.21}
	\begin{aligned}
		& \bigg|\frac{1}{2}\int_{\mathbb{R}^d}\left(\left|u(t, x)\right|^{2}\varphi(x)+\left|{\bf F}(t, x)\right|^{2}\phi(x)\right) d x-\frac{1}{2}\int_{\mathbb{R}^d}\left(\left|u(\tau, x)\right|^{2}\phi(x)+\left|{\bf F}(\tau, x)\right|^{2}\phi(x)\right) dx\\
		&+\int_{\tau}^{t}  \int_{\mathbb{R}^d}\left|\nabla u(s,x)\right|^{2} \phi(x)d x ds\bigg|\\
		=&\bigg|\frac{1}{2}\int_{\tau}^{t} \int_{\mathbb{R}^d}|u|^2u\cdot\nabla\phi dxds+ \int_{\tau}^{t} \int_{\mathbb{R}^d}\left\{P  \left(u\cdot \nabla\phi\right)-\nabla u:\left(u \otimes \nabla \phi\right)\right\}dxds\\
		&+\frac{1}{2}\int_{
			\tau}^{t} \int_{\mathbb{R}^d}u_k\partial_k\phi |{\bf F}_{ij}|^2dxds-\int_{\tau}^{t}\int_{\mathbb{R}^d}u_i{\bf F}_{ij}{\bf F}_{kj}\partial_k\phi dxds\bigg|\\
		=:& |\mathcal{F}_1+\mathcal{F}_2+\mathcal{F}_3+\mathcal{F}_4|.
	\end{aligned}
\end{equation}
Simply note that $$u \in L^{\infty}\left(0, T ; L^{2}\left(\mathbb{R}^d\right)\right) \cap L^{2}\left(0, T ; H^{1}\left(\mathbb{R}^d\right)\right),$$ from which it follows $$u \in L^{r}\left(0, T ; L^{s}\left(\mathbb{R}^d\right)\right),\,\, \frac{2}{r}+\frac{d}{s}=\frac{d}{2},\,\, 2\leq s< \frac{2d}{d-2}, \,\, d=2, 3, 4.$$
In particular, we have $$u \in L^{\frac{2d+4}{d}}\left(0, T ; L^{\frac{2d+4}{d}}\left(\mathbb{R}^d\right)\right)$$ and $$u \in L^{\frac{8d+16}{d^2}}\left(0, T ; L^{\frac{4d+8}{d+4}}\left(\mathbb{R}^d\right)\right)\subset L^{\frac{4d+8}{d+4}}\left(0, T ; L^{\frac{4d+8}{d+4}}\left(\mathbb{R}^d\right)\right).$$ Then by operating div on both
sides of the first equation in \eqref{1.1}, it follows that
$$
-\Delta P=\operatorname{div}\operatorname{div}\left(u \otimes u\right)-\operatorname{div}\operatorname{div} \left({\bf F }{\bf F}^{\top}\right)= \partial_{i} \partial_{k}\left(u_{i} u_{k}\right)+\partial_{i} \partial_{k}\left({\bf F }_{ij} {\bf F }_{kj}\right),
$$
from which together with the Calderon-Zygmund inequality yields for $1<s<\infty$
$$
\|P\|_{L^{s}} \leq C\left(\left\| u \otimes u\right\|_{L^{s}}+\left\|{\bf F }{\bf F}^{\top}\right\|_{L^{s}} \right)\leq C\left(\left\| u \right\|^2_{L^{2s}}+\left\|{\bf F }\right\|^2_{L^{2s}} \right).
$$
Since $(u, {\bf F } )\in L^{\frac{4d+8}{d+4}}\left(0, T ; L^{\frac{4d+8}{d+4}}\left(\mathbb{R}^d\right)\right)$,
i.e., $P \in L^{\frac{2d+4}{d+4}}\left(0, T ; L^{\frac{2d+4}{d+4}}\left(\mathbb{R}^d\right)\right)$.
Next, by using the H\"older's inequality, we can control the term related to $\mathcal{F}_1, \mathcal{F}_2, \mathcal{F}_3, \mathcal{F}_4$ in the following way
$$
\begin{aligned}
	&|\mathcal{F}_1+\mathcal{F}_2+\mathcal{F}_3+\mathcal{F}_4|\\
	\leq	& \frac{C}{r}\left(\|u\|_{L^{3}\left(0, T ; L^{3}\left(\mathbb{R}^d\right)\right)}^{3}+\|P\|_{L^{\frac{2d+4}{d+4}}\left(0, T ; L^{\frac{2d+4}{d+4}}\left(\mathbb{R}^d\right)\right)}\|u\|_{L^{\frac{2d+4}{d}}\left(0, T ; L^{\frac{2d+4}{d}}\left(\mathbb{R}^d\right)\right)}\right)\\
	+&\frac{C}{r}\left(\|u\|_{L^{\infty}\left(0, T ; L^{2}\left(\mathbb{R}^d\right)\right)}\|\nabla u\|_{L^{2}\left(0, T ; L^{2}\left(\mathbb{R}^d\right)\right)}+\|u\|_{L^{\frac{2d+4}{d}}\left(0, T ; L^{\frac{2d+4}{d}}\left(\mathbb{R}^d\right)\right)}\|{\bf F}\|^2_{L^{\frac{4d+8}{d+4}}\left(0, T ; L^{\frac{4d+8}{d+4}}\left(\mathbb{R}^d\right)\right)}\right).
\end{aligned}
$$
Then, by letting $r \rightarrow \infty$, one has
$$
|\mathcal{F}_1+\mathcal{F}_2+\mathcal{F}_3+\mathcal{F}_4| \rightarrow 0.
$$
Finally, letting $r \rightarrow \infty$ in \eqref{3.21}, and combining above convergence result, we obtain
$$
\frac{1}{2}\int_{\mathbb{R}^d}\left(\left|u(t, x)\right|^{2}+\left|{\bf F}(t, x)\right|^{2}\right) d x-\frac{1}{2}\int_{\mathbb{R}^d}\left(\left|u(\tau, x)\right|^{2}+\left|{\bf F}(\tau, x)\right|^{2}\right) dx+ \int_{\tau}^{t}  \int_{\mathbb{R}^d}\left|\nabla u(s,x)\right|^{2} d x ds=0
$$
for $0<\tau\leq t\leq T$. By passing to the limit as $\tau\rightarrow 0$ and using the facts that $\lim_{t\to 0^+}\|u(\cdot,t)\|^2_{L^2}+\|{\bf F}(\cdot,t)\|^2_{L^2}=\|u_0\|^2_{L^2}+\|{\bf F}_0\|^2_{L^2}$,  will yield energy equality \eqref{1.2}.
We conclude the proof of Theorem \ref{th2}.
\section{On the energy equality for distributional solutions.}
In this section, we shall give the proof of Theorem \ref{th3}. To this end, we need to introduce
a crucial proposition.
\begin{proposition}\label{pro4.1}
	Let $(u_0, {\bf F}_0)\in L^2(\mathbb{R}^{d})$ with $\nabla\cdot u_0 = 0$ and let $(u, {\bf F})$ be a distributional solution in the sense of Definition \ref{de2} to system \eqref{1.1} and satisfies
	$$(u, {\bf F}) \in L^{4}\left(0, T ; L^{4}(\mathbb{R}^{d})\right),$$ then we have
	$$
	\sup _{t \geq 0}\|u^\varepsilon(\cdot, t)\|_{L^{2}}^{2}+\frac{3}{2}\int^t_0\int_{\mathbb{R}^{d}}|\nabla u^\varepsilon|^{2} dx d\tau\leq C,\quad \forall\, t \in[0, T],
	$$
	where $C$ is a constant depending only on $\|u_0\|_{L^2}$ and $\int^T_0\|u\|_{L^{4}}^{4}+\|{\bf F}\|_{L^{4}}^{4} dt$.
\end{proposition}
\begin{remark}
Proposition \ref{pro4.1} shows that the distributional solution $u$ can be identified with $u \in L^{\infty}\left(0, T ; L^{2}(\mathbb{R}^d)\right) \cap L^{2}\left(0, T ; H^{1}(\mathbb{R}^d)\right)$. This proves that $u$ falls into the class of Leray-Hopf weak solutions, for which the classical results from \cite{SHIN} imply the energy equality to Navier-Stokes equations.  We note that the Galdi's result \cite{GALDI} can be treated directly by this proposition.
\end{remark}
\textbf{Proof of Proposition \ref{pro4.1}}. By the definition of distributional solutions, we obtain that following identity
$$
\begin{aligned}
	&\int_{\mathbb{R}^{d}} u \cdot \partial_{t} \Phi^\varepsilon+u \cdot \Delta \Phi^\varepsilon+u \otimes u: \nabla \otimes \Phi^\varepsilon-{\bf F}{\bf F}^{\top}: \nabla \otimes \Phi^\varepsilon d x=\frac{d}{dt}\int_{\mathbb{R}^{d}} u(x, t) \cdot \Phi^\varepsilon(x, t) d x,
\end{aligned}
$$
for all $\Phi^\varepsilon\in \mathcal{D}_{T}$. Which in turn gives
$$
\begin{aligned}
	&\int_{\mathbb{R}^{d}} u^\varepsilon \cdot \partial_{t} \Phi+u^\varepsilon \cdot \Delta \Phi+(u \otimes u)^\varepsilon: \nabla \otimes \Phi-({\bf F}{\bf F}^{\top})^\varepsilon: \nabla \otimes \Phi d x dt=\frac{d}{dt}\int_{\mathbb{R}^{d}} u^\varepsilon(x, t) \cdot \Phi(x, t) d x.
\end{aligned}
$$
Now, choosing $\Phi=u^\varepsilon$ in above identity, integrate by parts to find
\begin{equation}
	\begin{split}\label{3.22}
		\frac{1}{2}\frac{d}{dt}\int_{\mathbb{R}^{d}}| u^\varepsilon|^{2} dx+\int_{\mathbb{R}^{d}}|\nabla u^\varepsilon|^{2} dx
		&=\int_{\mathbb{R}^{d}}(u\otimes u)^\varepsilon\cdot\nabla u^\varepsilon dx-\int_{\mathbb{R}^{d}}({\bf F}{\bf F}^{\top})^\varepsilon\cdot \nabla u^\varepsilon dx\\
		&=I_1+I_2.
	\end{split}
\end{equation}
For $I_1$, applying the H\"older inequality, we have
\begin{equation}
	\begin{split}\label{3.23}
		I_{1} &\leq\left|\int_{\mathbb{R}^{d}}(u\otimes u)^\varepsilon\cdot\nabla u^\varepsilon dx\right|\leq C\|(u\otimes u)^\varepsilon\|_{L^2}\|\nabla u^\varepsilon\|_{L^{2}}\\
		& \leq C\|(u\otimes u)\|_{L^2}\|\nabla u^\varepsilon\|_{L^{2}}\leq C\|u\|_{L^{4}}^{2}\|\nabla u^\varepsilon\|_{L^{2}}\\
		&\leq C\|u\|_{L^{4}}^{4}+\frac{1}{4}\|\nabla u^\varepsilon\|^2_{L^{2}}.
	\end{split}
\end{equation}
Similarily, for $I_2$, we obtain
\begin{equation}
	\begin{split}\label{3.24}
		I_{2} &\leq\left|\int_{\mathbb{R}^{d}}({\bf F}{\bf F}^{\top})^\varepsilon\cdot\nabla u^\varepsilon dx\right| \leq C\|{\bf F}\|_{L^{4}}^{2}\|\nabla u^\varepsilon\|_{L^{2}}\\
		&\leq C\|{\bf F }\|_{L^{4}}^{4}+\frac{1}{4}\|\nabla u^\varepsilon\|^2_{L^{2}}.
	\end{split}
\end{equation}
Putting the above estimates \eqref{3.23}-\eqref{3.24} into \eqref{3.22}, one concludes that
\begin{equation}
	\begin{split}\label{3.25}
		\frac{d}{dt}\int_{\mathbb{R}^{d}}| u^\varepsilon|^{2} dx+\frac{3}{2}\int_{\mathbb{R}^{d}}|\nabla u^\varepsilon|^{2} dx\leq C\left(\|u\|_{L^{4}}^{4}+\|{\bf F}\|_{L^{4}}^{4}\right),
	\end{split}
\end{equation}
and it follows that
\begin{equation}\label{3.26}
	\sup _{t \geq 0}\|u^\varepsilon(\cdot, t)\|_{L^{2}}^{2}+\frac{3}{2}\int^t_0\int_{\mathbb{R}^{d}}|\nabla u^\varepsilon|^{2} dx d\tau \leq\left\|u_{0}\right\|_{L^{2}}^{2}+C\int^t_0
	\left(\|u\|_{L^{4}}^{4}+\|{\bf F}\|_{L^{4}}^{4}\right) d\tau\leq C
\end{equation}
for all $t\in [0, T]$, where $C$ is a constant depending only on $u_0$ and $\int^T_0\|u\|_{L^{4}}^{4}+\|{\bf F}\|_{L^{4}}^{4} dt$.
Moreover, by the $L^p$-norm is weakly lower semicontinuous, one has
\begin{equation}\label{3.27}
	\sup _{t \geq 0}\|u(\cdot, t)\|_{L^{2}}^{2}+\frac{3}{2}\int^t_0\int_{\mathbb{R}^{d}}|\nabla u|^{2} dx d\tau\leq C.
\end{equation}
This concludes the proof of Proposition \ref{pro4.1}.

\textbf{Proof of Theorem \ref{th3}}.
With  Proposition \ref{pro4.1}  in hand, we are ready to prove Theorem \ref{th3}. First, we define two new functions $\Xi=u^\varepsilon$,  $\Sigma={\bf F}^\varepsilon$,  and note that, we have
$$\operatorname{div} u^\varepsilon=0. $$
Using $\Xi^\varepsilon$ and $\Sigma^\varepsilon$ to test $\eqref{1.1}_1$ and $\eqref{1.1}_2$, respectively, one has
\begin{equation}\label{3.28}
	\left\{\begin{array}{lr}
		\int_{\mathbb{R}^{d}}	\Xi^\varepsilon\left(\partial_{t} u+u \cdot \nabla u- \Delta u+\nabla P-\nabla\cdot({\bf F}{\bf F}^\top)\right) dx=0, \\
		\int_{\mathbb{R}^{d}}\Sigma^\varepsilon\left( \partial_t {\bf F}+(u\cdot\nabla) {\bf F}-\nabla u {\bf F}\right) dx=0,
	\end{array}\right.
\end{equation}
we get that
\begin{equation}\label{3.29}
	\left\{\begin{array}{lr}
		\int_{\mathbb{R}^{d}}	u^\varepsilon\left(\partial_{t} u+u \cdot \nabla u- \Delta u+\nabla P-\nabla\cdot({\bf F}{\bf F}^\top)\right)^\varepsilon dx=0, \\
		\int_{\mathbb{R}^{d}}{\bf F}^\varepsilon\left(\partial_t {\bf F}+(u\cdot\nabla) {\bf F}-\nabla u {\bf F}\right)^\varepsilon dx=0.
	\end{array}\right.
\end{equation}
This yields
\begin{equation}
	\begin{split}\label{3.30}
		&\frac{1}{2}\frac{d}{dt}\int_{\mathbb{R}^{d}}\left(| u^\varepsilon|^{2}+|{\bf F}^\varepsilon|^{2}\right) dx+\int_{\mathbb{R}^{d}}|\nabla u^\varepsilon|^{2} dx\\
		=&-\int_{\mathbb{R}^{d}}\operatorname{div}(u\otimes u)^\varepsilon\cdot u^\varepsilon dx+\int_{\mathbb{R}^{d}}\operatorname{div}({\bf F}{\bf F}^\top)^\varepsilon\cdot u^\varepsilon dx-\int_{\mathbb{R}^{d}}(u\cdot\nabla {\bf F})^\varepsilon\cdot {\bf F}^\varepsilon dx\\
		&+\int_{\mathbb{R}^{d}}( \nabla u{\bf F})^\varepsilon \cdot {\bf F}^\varepsilon dx.
	\end{split}
\end{equation}
Clearly,
\begin{equation}
	\begin{split}\label{3.31}
		&\int_{\mathbb{R}^{d}}\left(| u^\varepsilon(x,t)|^{2}+|{\bf F}^\varepsilon(x,t)|^{2}\right) dx-\int_{\mathbb{R}^{d}}\left(| u_0^\varepsilon(x,0)|^{2}+|{\bf F}_0^\varepsilon(x,0)|^{2}\right) dx+2\int^t_0\int_{\mathbb{R}^{d}}|\nabla u^\varepsilon|^{2} dxd\tau\\
		=&-2\int^t_0\int_{\mathbb{R}^{d}}\operatorname{div}(u\otimes u)^\varepsilon\cdot u^\varepsilon dxd\tau+2\int^t_0\int_{\mathbb{R}^{d}}\operatorname{div}({\bf F}{\bf F}^\top)^\varepsilon\cdot u^\varepsilon dxd\tau-2\int^t_0\int_{\mathbb{R}^{d}}(u\cdot\nabla {\bf F})^\varepsilon\cdot {\bf F}^\varepsilon dxd\tau\\
		&+2\int^t_0\int_{\mathbb{R}^{d}}( \nabla u{\bf F})^\varepsilon \cdot {\bf F}^\varepsilon dxd\tau\\
		=&\mathcal{J}_1+\mathcal{J}_2+\mathcal{J}_3+\mathcal{J}_4.
	\end{split}
\end{equation}
Notice that
$$
-2\int^t_0\int_{\mathbb{R}^{d}}\operatorname{div}(u^\varepsilon \otimes u^\varepsilon)\cdot u^\varepsilon dx d\tau=0,
$$
thus by using the H\"oder's equality, one has
\begin{equation}
	\begin{split}\label{3.32}
		\mathcal{J}_1=&-2\int^t_0\int_{\mathbb{R}^{d}}\operatorname{div}(u\otimes u)^\varepsilon\cdot u^\varepsilon-\operatorname{div}(u^\varepsilon \otimes u^\varepsilon)\cdot u^\varepsilon dxd\tau\\
		=&2\int^t_0\int_{\mathbb{R}^{d}}[(u\otimes u)^\varepsilon-(u^\varepsilon \otimes u^\varepsilon)]\cdot\nabla u^\varepsilon dxd\tau\\
		\leq &2\int_{0}^{t} \int_{\mathbb{R}^{d}}\left|(u \otimes u)^\varepsilon-u^\varepsilon \otimes u^\varepsilon\right| \left| \nabla u^\varepsilon\right| d x d \tau\\
		\leq &2\int_{0}^{t} \int_{\mathbb{R}^{d}}\left(\left|(u \otimes u)^\varepsilon-u \otimes u\right|+\left|u \otimes u-u \otimes u^\varepsilon\right|+\left|u \otimes u^\varepsilon-u^\varepsilon \otimes u^\varepsilon\right|\right) \left| \nabla u^\varepsilon\right| d x d \tau\\
		\leq&C\left\|(u \otimes u)^\varepsilon-u \otimes u\right\|_{L^{2}\left(0, T ; L^{2}(\mathbb{R}^{d})\right)} \|\nabla u^\varepsilon \|_{L^{2}\left(0, T ; L^{2}(\mathbb{R}^{d})\right)}\\
		&+C\left\|u-u^\varepsilon\right\|_{L^{4}\left(0, T ; L^{4}(\mathbb{R}^{d})\right)} \|u \|_{L^{4}\left(0, T ; L^{4}(\mathbb{R}^{d})\right)} \|\nabla u^\varepsilon \|_{L^{2}\left(0, T ; L^{2}(\mathbb{R}^{d})\right)}\\
		&+C\left\|u-u^\varepsilon\right\|_{L^{4}\left(0, T ; L^{4}(\mathbb{R}^{d})\right)} \|u^\varepsilon \|_{L^{4}\left(0, T ; L^{4}(\mathbb{R}^{d})\right)} \|\nabla u^\varepsilon \|_{L^{2}\left(0, T ; L^{2}(\mathbb{R}^{d})\right)}\\
		&	\rightarrow  0, \quad as\quad \epsilon\rightarrow 0.
	\end{split}
\end{equation}
By using H\"older's inequality and Lemma \ref{le2.3}, we obtain
\begin{equation}\label{3.33}
	\begin{split}
		\mathcal{J}_3&=-2\int^t_0\int_{\mathbb{R}^{d}}(u\cdot\nabla {\bf F})^\varepsilon\cdot {\bf F}^\varepsilon dxd\tau\\
		&=-2\int^t_0\int_{\mathbb{R}^{d}}\partial_k(u_k{\bf F}_{ij})^{\varepsilon}\cdot {\bf F}_{ij}^\varepsilon dxd\tau\\
		&=2\int^t_0\int_{\mathbb{R}^{d}}(u_k{\bf F}_{ij})^{\varepsilon}\partial_k {\bf F}_{ij}^\varepsilon dxd\tau\\
		&=2\int^t_0\int_{\mathbb{R}^{d}}[(u_k{\bf F}_{ij})^{\varepsilon}-u^{\varepsilon}_k{\bf F}^{\varepsilon}_{ij}]\partial_k {\bf F}_{ij}^\varepsilon dxd\tau\\
		&\leq C\varepsilon^{-1}\left\|(u_k{\bf F}_{ij})^{\varepsilon}-u^{\varepsilon}_k{\bf F}^{\varepsilon}_{ij}\right\|_{L^{\frac{4}{3}}\left(0, T ; L^{\frac{4}{3}}(\mathbb{R}^{d})\right)} \|{\bf F}_{ij}^\varepsilon \|_{L^{4}\left(0, T ; L^{4}(\mathbb{R}^{d})\right)} \\
		&	\rightarrow  0, \quad as\quad \varepsilon\rightarrow 0,
	\end{split}
\end{equation}
where we have used the fact
$$\int^t_0\int_{\mathbb{R}^{d}}u^{\varepsilon}_k{\bf F}^{\varepsilon}_{ij}\partial_k {\bf F}_{ij}^\varepsilon dxd\tau=0.$$
Now, we need to show that, when $\varepsilon\rightarrow 0$, the term $\mathcal{J}_2$ and $\mathcal{J}_4$ have the following limit:

\begin{equation}\label{3.34}
	\begin{split}
		\mathcal{J}_2&=2\int^t_0\int_{\mathbb{R}^{d}}\operatorname{div}({\bf F}{\bf F}^{\top})^\varepsilon\cdot u^\varepsilon dxdt=-2\int^t_0\int_{\mathbb{R}^{d}}({\bf F}{\bf F}^{\top})^\varepsilon\cdot \nabla u^\varepsilon dxd\tau\\
		&=-2\int^t_0\int_{\mathbb{R}^{d}}({\bf F}_{ij}{\bf F}_{kj})^\varepsilon\partial_k u_i^\varepsilon dxd\tau\\
			&	\rightarrow  -2\int^t_0\int_{\mathbb{R}^{d}}{\bf F}_{ij}{\bf F}_{kj}\partial_k u_i dxd\tau, \quad as\quad \varepsilon\rightarrow 0
	\end{split}
\end{equation}
and
\begin{equation}\label{3.35}
	\begin{split}
		\mathcal{J}_4&=2\int^t_0\int_{\mathbb{R}^{d}}( \nabla u{\bf F})^\varepsilon \cdot {\bf F}^\varepsilon dxd\tau\\
		&=2\int^t_0\int_{\mathbb{R}^{d}}(\partial_k u_i{\bf F}_{kj})^\varepsilon ({\bf F}_{ij})^\varepsilon dxd\tau\\
		&	\rightarrow  2\int^t_0\int_{\mathbb{R}^{d}}\partial_k u_i{\bf F}_{kj}{\bf F}_{ij} dxd\tau, \quad as\quad \varepsilon\rightarrow 0.
	\end{split}
\end{equation}
Indeed, from the the uniform bound of $u^\varepsilon$ in $L^\infty L^2\cap L^2H^1$ established in Proposition \ref{pro4.1}, we have
\begin{equation}\label{3.36}
	\begin{split}
&\left|2\int^t_0\int_{\mathbb{R}^{d}}({\bf F}_{ij}{\bf F}_{kj})^\varepsilon\partial_k u_i^\varepsilon-{\bf F}_{ij}{\bf F}_{kj}\partial_k u_i dxd\tau\right| \\
= & \left|2\int^t_0\int_{\mathbb{R}^{d}}[({\bf F}_{ij}{\bf F}_{kj})^{\varepsilon}-{\bf F}_{ij}{\bf F}_{kj}]\partial_ku^{\varepsilon}_i+[\partial_ku^{\varepsilon}_i-\partial_ku_i]{\bf F}_{ij}{\bf F}_{kj} dx d\tau\right|\\
\leq &2 \int^t_0\int_{\mathbb{R}^{d}}\left|({\bf F}_{ij}{\bf F}_{kj})^{\varepsilon}-{\bf F}_{ij}{\bf F}_{kj}\right|\left|\partial_ku^{\varepsilon}_i\right| dx d\tau+2 \int^t_0\int_{\mathbb{R}^{d}}\left|\partial_ku^{\varepsilon}_i-\partial_ku_i\right|\left|{\bf F}_{ij}{\bf F}_{kj}\right| dx d\tau\\
\leq&C\left\|({\bf F}_{ij}{\bf F}_{kj})^{\varepsilon}-{\bf F}_{ij}{\bf F}_{kj}\right\|_{L^{2}\left(0, T ; L^{2}(\mathbb{R}^{d})\right)} \|\partial_ku^{\varepsilon}_i \|_{L^{2}\left(0, T ; L^{2}(\mathbb{R}^{d})\right)}\\
&+C\left\|\partial_ku^{\varepsilon}_i-\partial_ku_i\right\|_{L^{2}\left(0, T ; L^{2}(\mathbb{R}^{d})\right)} \|{\bf F}_{ij}{\bf F}_{kj} \|_{L^{2}\left(0, T ; L^{2}(\mathbb{R}^{d})\right)}\\
&\rightarrow 0,\,\, as \,\,\varepsilon \rightarrow 0
\end{split}
\end{equation}
and
\begin{equation}\label{3.37}
	\begin{split}
		&\left|2\int^t_0\int_{\mathbb{R}^{d}}(\partial_k u_i{\bf F}_{kj})^\varepsilon ({\bf F}_{ij})^\varepsilon-\partial_k u_i{\bf F}_{kj} {\bf F}_{ij} dxd\tau\right| \\
		= & \left|2\int^t_0\int_{\mathbb{R}^{d}}[{\bf F}_{ij}^{\varepsilon}-{\bf F}_{ij}](\partial_ku_i{\bf F}_{kj})^{\varepsilon}+[(\partial_ku_i{\bf F}_{kj})^{\varepsilon}-\partial_ku_i{\bf F}_{kj}]{\bf F}_{ij} dx d\tau\right|\\
		\leq&C\left\|{\bf F}_{ij}^{\varepsilon}-{\bf F}_{ij}\right\|_{L^{4}\left(0, T ; L^{4}(\mathbb{R}^{d})\right)} \|(\partial_ku_i{\bf F}_{kj})^{\varepsilon} \|_{L^{\frac{4}{3}}\left(0, T ; L^{\frac{4}{3}}(\mathbb{R}^{d})\right)}\\
		&+C\left\|(\partial_ku_i{\bf F}_{kj})^{\varepsilon}-\partial_ku_i{\bf F}_{kj}\right\|_{L^{\frac{4}{3}}\left(0, T ; L^{\frac{4}{3}}(\mathbb{R}^{d})\right)} \|{\bf F}_{ij} \|_{L^{4}\left(0, T ; L^{4}(\mathbb{R}^{d})\right)}\\
		&\rightarrow 0,\,\, as \,\,\varepsilon \rightarrow 0,
	\end{split}
\end{equation}
which completes the proof of \eqref{3.34}-\eqref{3.35}.

Then it follows from \eqref{3.31}-\eqref{3.35} that
\begin{equation}\label{3.38}
\left|\mathcal{J}_1+	\mathcal{J}_2+\mathcal{J}_3+\mathcal{J}_4\right|\rightarrow 0,  \quad as\quad \epsilon\rightarrow 0.
\end{equation}
Letting $\varepsilon$ goes to zero in \eqref{3.31}, and using the facts \eqref{3.38} yields
\begin{equation}
	\begin{split}\label{3.39}
		\int_{\mathbb{R}^{d}}\left(| u(x,t)|^{2}+|{\bf F}(x,t)|^{2}\right) dx+2\int^t_0\int_{\mathbb{R}^{d}}|\nabla u|^{2} dxd\tau=\int_{\mathbb{R}^{d}}\left(| u_0|^{2}+|{\bf F}_0|^{2}\right) dx,
	\end{split}
\end{equation}
so we obtain the assertion of Theorem \ref{th3}.

\section*{Acknowledgments}

The authors are supported by the Construct Program of the Key Discipline in Hunan Province and NSFC Grant No. 11871209.

\section*{Data Availability}
The data that support the findings of this study are available from the corresponding author upon reasonable request.


\begin{thebibliography}{10}
\bibitem{AH}
Amann H. On the strong solvability of the Navier-Stokes equations. Journal of Mathematical Fluid Mechanics, 2000, 2(1): 16-98.
\bibitem{BT}
Bardos C, Titi E. Onsager's conjecture for the incompressible Euler equations in bounded domains. Archive for Rational Mechanics and Analysis, 2018, 228(1): 197-207.
\bibitem{BLC}
Berselli L, Chiodaroli E. On the energy equality for the 3D Navier-Stokes equations. Nonlinear Analysis, 2020, 192: 111704.
\bibitem{BG}
Berselli L, Galdi G. On the space–time regularity of $C(0,T;L^n)$-very weak solutions to the Navier-Stokes equations. Nonlinear Analysis: Theory, Methods \& Applications, 2004, 58(5-6): 703-717.
\bibitem{Vicol}
Buckmaster T, Vicol V. Non-uniqueness of weak solutions to the Navier-Stokes equations. Annals of Mathematics, 2019, 189 101-144.
\bibitem{CPT}
Constantin P, E W, Titi E. Onsager's conjecture on the energy conservation for solutions of Euler's equation. Communications in Mathematical Physics, 1994, 165(1): 207-209.
\bibitem{CLWX}
Chen M, Liang Z, Wang D, Xu R. Energy equality in compressible fluids with physical boundaries. SIAM Journal on Mathematical Analysis, 2020, 52(2): 1363-1385.
\bibitem{CW}
Chen Q, Wu G. The 3D compressible viscoelastic fluid in a bounded domain. Communications in Mathematical Sciences, 2018, 16(5): 1303-1323.
\bibitem{CCFS}
Cheskidov A, Constantin P, Friedlander S, Shvydkoy R. Energy conservation and Onsager's conjecture for the Euler equations. Nonlinearity, 2008, 21(6): 1233-1252.
\bibitem{CF}
Cheskidov A, Friedlander S, Shvydkoy R. On the energy equality for weak solutions of the 3D Navier-Stokes equations. Advances in mathematical fluid mechanics. Springer, Berlin, Heidelberg, 2010: 171-175.
\bibitem{CLX}
Cheskidov A, Luo X. Energy equality for the Navier-Stokes equations in weak-in-time Onsager spaces. Nonlinearity, 2020, 33(4): 1388-1403.

\bibitem{EGL}
Eyink G L. Energy dissipation without viscosity in ideal hydrodynamics I. Fourier analysis and local energy transfer. Physica D: Nonlinear Phenomena, 1994, 78(3-4): 222-240.
\bibitem{ELSGA}
Escauriaza L, Seregin G, Sverak V. $L_{3,\infty}$-solutions of the Navier-Stokes equations and backward uniqueness. Russian Mathematical Surveys, 2003, 58(2): 211-250.
\bibitem{FZZ}
Feng Z, Zhu C, Zi R. Blow-up criterion for the incompressible viscoelastic flows. Journal of Functional Analysis, 2017, 272(9): 3742-3762.
\bibitem{FJR}
Fabes E, Jones B, Riviere, N. The initial value problem for the Navier-Stokes equations with data in $L^p$. Archive for Rational Mechanics and Analysis, 1972, 45, 222-240.
\bibitem{GALDI1}
Galdi G. An introduction to the mathematical theory of the Navier-Stokes equations: Steady-state problems. Springer Science \& Business Media, 2011.
\bibitem{GALDI}
Galdi G. On the energy equality for distributional solutions to Navier-Stokes equations. Proceedings of the American Mathematical Society, 2019, 147(2): 785-792.
\bibitem{GALDI2}
Galdi G. On the relation between very weak and Leray–Hopf solutions to Navier-Stokes equations. Proceedings of the American Mathematical Society, 2019, 147(12): 5349-5359.
\bibitem{HZ}
He Y, Zi R. Energy conservation for solutions of incompressible viscoelastic fluids. Acta Mathematica Scientia, 2021, 41(4): 1287-1301.
\bibitem{HXP}
Hu X. Global existence of weak solutions to two dimensional compressible viscoelastic flows. Journal of Differential Equations, 2018, 265(7): 3130-3167.
\bibitem{HL}
Hu X, Lin F. On the Cauchy problem for two dimensional incompressible viscoelastic flows. arXiv preprint arXiv:1601.03497, 2016.
\bibitem{HUW}
Hu X, Wu G. Global existence and optimal decay rates for three-dimensional compressible viscoelastic flows. SIAM Journal on Mathematical Analysis, 2013, 45(5): 2815-2833.
\bibitem{HW}
Hu X, Wu H. Long-time behavior and weak-strong uniqueness for incompressible viscoelastic flows. Discrete and Continuous Dynamical Systems 2015,35(8): 3437-3461.

\bibitem{HH}
Hu X, Hynd R. A blowup criterion for ideal viscoelastic flow. Journal of Mathematical Fluid Mechanics 2013,15(3): 431-437.
\bibitem{IY}
Ishigaki Y. Diffusion wave phenomena and $L^p$ decay estimates of solutions of compressible viscoelastic system. Journal of Differential Equations, 2020, 269(12): 11195-11230.
\bibitem{IS}
Isett P. A proof of Onsager's conjecture. Annals of Mathematics, 2018, 188(3): 871-963.

\bibitem{FC}
Foias C. Une remarque sur l'unicit\'e des solutions des \'equations de
{N}avier-{S}tokes en dimension {$n$}. Bulletin de la Soci\'et\'e Math\'ematique de France, 1961, 89: 1-8.
\bibitem{KT}
Kozono H, Taniuchi Y. Bilinear estimates in BMO and the Navier-Stokes equations. Math Z, 2000, 235: 173-194.

\bibitem{RGL}
Larson, R. The Structure and Rheology of Complex Fluids, Oxford University Press, New York, 1995.

\bibitem{LZL}
Lei Z, Liu C, Zhou Y. Global solutions for incompressible viscoelastic fluids. Archive for Rational Mechanics and Analysis, 2008, 188(3): 371-398.
\bibitem{LS}
Leslie T, Shvydkoy R. The energy measure for the Euler and Navier-Stokes equations. Archive for Rational Mechanics and Analysis, 2018, 230(2): 459-492.
\bibitem{LS1}
Leslie T, Shvydkoy R. Conditions implying energy equality for weak solutions of the Navier-Stokes equations. SIAM Journal on Mathematical Analysis, 2018, 50(1): 870-890.
\bibitem{LLZ}
Lin F, Liu C, Zhang P. On hydrodynamics of viscoelastic fluids. Communications on Pure and Applied Mathematics, 2005, 58(11): 1437-1471.
\bibitem{LINFANGHUA}
Lin F. Some analytical issues for elastic complex fluids. Communications on Pure and Applied Mathematics, 2012, 65(7): 893-919.
\bibitem{LZ}
Lin F, Zhang P. Global small solutions to an MHD-type system: the three-dimensional case. Communications on Pure and Applied Mathematics, 2014, 67(4): 531-580.
\bibitem{LXZ}
Lin F, Xu L, Zhang P. Global small solutions of 2-D incompressible MHD system. Journal of Differential Equations, 2015, 259(10): 5440-5485.

\bibitem{LIONS}
Lions J. Sur la r\'egularit\'e et l'unicit\'e des solutions turbulentes des \'equations de Navier-Stokes. Rendiconti del Seminario Matematico della Universita di Padova, 1960, 30: 16-23.
\bibitem{LM}
Lions P, Masmoudi N. Uniqueness of mild solutions of the Navier-Stokes system in $L^N$. Communications in Partial Differential Equations, 2001, 26(11-12): 2211-2226.
\bibitem{NQH}
Nguyen Q, Nguyen P, Tang B. Energy equalities for compressible Navier-Stokes equations. Nonlinearity, 2019, 32(11): 4206-4231.
\bibitem{ON}
Onsager L. Statistical Hydrodynamics, Nuovo Cimento (Supplemento) 6 (1949), 279-287.
\bibitem{GP}
Prodi G. Un teorema di unicita per le equazioni di Navier-Stokes. Annali di Matematica pura ed applicata, 1959, 48(1): 173-182.
\bibitem{RM}
Renardy M. Mathematical analysis of viscoelastic flows. Society for Industrial and Applied Mathematics, 2000.
\bibitem{SS}
Seregin G, Sverak V. Navier-Stokes equations with lower bounds on the pressure.  Archive for Rational Mechanics and Analysis, 2002, 163(1): 65-86.
\bibitem{SJ}
Serrin J. On the interior regularity of weak solutions of the Navier-Stokes equations. Archive for Rational Mechanics and Analysis, 1962, 9(1): 187-195.
\bibitem{SHIN}
Shinbrot M. The energy equation for the Navier-Stokes system. SIAM Journal on Mathematical Analysis, 1974, 5(6): 948-954.

\bibitem{SOHR}
Sohr H. The Navier-Stokes equations: An elementary functional analytic approach. Springer Science \& Business Media, 2012.
\bibitem{TAN}
Tan W, Yin Z. The energy conservation and regularity for the Navier-Stokes equations. arXiv preprint arXiv:2107.04157, 2021.
\bibitem{YUCHENG}
Yu C. A new proof of the energy conservation for the Navier-Stokes equations. arxiv:1604.05697.
\bibitem{YC}
Yu C. The energy equality for the Navier-Stokes equations in bounded domains. arXiv preprint arXiv:1802.07661, 2018.
\bibitem{YBQ}
Yuan B. Note on the blowup criterion of smooth solution to the incompressible viscoelastic flow. Discrete and Continuous Dynamical Systems 2013,33(5):2211-2219.


\end{thebibliography}
\end{document}